
\magnification=1200
\def\'#1{{\if #1i{\accent"13\i}\else {\accent"13 #1}\fi}}

\hsize=32pc
\vsize=43pc
\vsize 19truecm
\baselineskip 14pt

\global\newcount\numsec\global\newcount\numfor

\def\senondefinito#1{\expandafter\ifx\csname#1\endcsname\relax}

\def\SIA #1,#2,#3 {\senondefinito{#1#2}%
\expandafter\xdef\csname #1#2\endcsname{#3}\else
\write16{???? ma #1,#2 e' gia' stato definito !!!!} \fi}

\def\etichetta(#1){(\veroparagrafo.\veraformula)%
\SIA e,#1,(\veroparagrafo.\veraformula) %
\global\advance\numfor by 1%
\write16{ EQ #1 ==> \equ(#1) }}

\def\letichetta(#1){\veroparagrafo.\veraformula
\SIA e,#1,{\veroparagrafo.\veraformula}
\global\advance\numfor by 1
 \write16{ Sta #1 ==> \equ(#1) }}

\def\tetichetta(#1){\veroparagrafo.\veraformula 
\SIA e,#1,{(\veroparagrafo.\veraformula)}
\global\advance\numfor by 1
 \write16{ tag #1 ==> \equ(#1)}}

\def\FU(#1)#2{\SIA fu,#1,#2 }

\def\etichettaa(#1){(A\veroparagrafo.\veraformula)%
\SIA e,#1,(A\veroparagrafo.\veraformula) %
\global\advance\numfor by 1%
\write16{ EQ #1 ==> \equ(#1) }}

\def\BOZZA{
\def\alato(##1){%
 {\rlap{\kern-\hsize\kern-1.4truecm{$\scriptstyle##1$}}}}%
\def\aolado(##1){%
 {
{
 \rlap{\kern-1.4truecm{$\scriptstyle##1$}}}}}
}

\def\alato(#1){}
\def\aolado(#1){}

\def\veroparagrafo{\number\numsec}\def\veraformula{\number\numfor}

\def\Eq(#1){\eqno{\etichetta(#1)\alato(#1)}}
\def\eq(#1){\etichetta(#1)\alato(#1)}
\def\leq(#1){\leqno{\aolado(#1)\etichetta(#1)}}
\def\teq(#1){\tag{\aolado(#1)\tetichetta(#1)\alato(#1)}}
\def\Eqa(#1){\eqno{\etichettaa(#1)\alato(#1)}}
\def\eqa(#1){\etichettaa(#1)\alato(#1)}
\def\eqv(#1){\senondefinito{fu#1}$\clubsuit$#1
\write16{#1 non e' (ancora) definito}%
\else\csname fu#1\endcsname\fi}
\def\equ(#1){\senondefinito{e#1}\eqv(#1)\else\csname e#1\endcsname\fi}

\def\Lemma(#1){\aolado(#1)Lemma \letichetta(#1)}%
\def\Theorem(#1){{\aolado(#1)Theorem \letichetta(#1)}}%
\def\Proposition(#1){\aolado(#1){Proposition \letichetta(#1)}}%
\def\Corollary(#1){{\aolado(#1)Corollary \letichetta(#1)}}%
\def\Remark(#1){\aolado(#1)Remark \letichetta(#1)}%

\def\include#1{
\openin13=#1.aux \ifeof13 \relax \else
\input #1.aux \closein13 \fi}

\openin14=\jobname.aux \ifeof14 \relax \else
\input \jobname.aux \closein14 \fi
\openout15=\jobname.aux

\def\include#1{
\openin13=#1.aux \ifeof13 \relax \else
\input #1.aux \closein13 \fi}
\openin14=\jobname.aux \ifeof14 \relax \else
\input \jobname.aux \closein14 \fi

\def\one{{\bf 1}}
\def\ni{{\noindent }}
\def\sobre#1#2{\lower 1ex \hbox{ $#1 \atop #2 $ } }
\input amssym.def  
\input amssym.tex  
\def\E{{\Bbb E}}

\def\P{{\Bbb P}}

\def\vep{{\varepsilon}}
\def\Z{{\Bbb Z}}

\def\N{{\Bbb N}}
\def\bn{{\bf N}}
\def\bbn{{\bf \overline N}}
\def\R{{\Bbb R}}
\def\bm{{\bf M}}
\def\({\left(}
\def\){\right)}
\def\square{\ifmmode\sqr\else{$\sqr$}\fi}
\def\sqr{\vcenter{
         \hrule height.1mm
         \hbox{\vrule width.1mm height2.2mm\kern2.18mm\vrule width.1mm}
         \hrule height.1mm}}                  

\bigskip
\centerline{\bf
 Ces\`aro 
mean distribution of group automata}
\centerline{\bf starting from measures with summable decay}
\vskip 3mm
\centerline{\bf Pablo A. Ferrari}
\centerline{\sl
Instituto de Matem\'atica e Estat\'istica,
Universidade de S\~ao Paulo}
\centerline{\sl
Caixa Postal 66281, 05315-970 S\~ao Paulo,
Brasil}
\centerline{\sl e-mail: \tt pablo@ime.usp.br}
\vskip 3mm
\centerline{\bf Alejandro Maass, Servet Mart\'inez}
\centerline{\sl
Departamento de Ingenier\'{\i}a Matem\'atica}
\centerline{\sl Universidad de Chile,
Facultad de Ciencias F\'{\i}sicas y
Matem\'aticas,}
\centerline{\sl
Casilla 170-3, Correo 3, Santiago, Chile}
\centerline{\sl e-mail: \tt 
amaass@dim.uchile.cl, smartine@dim.uchile.cl}
\vskip 3mm
\centerline{\bf Peter Ney}
\centerline{\sl
Department of Mathematics, University of Wisconsin,}
\centerline{\sl
Madison, WI 53706, U.S.A.}
\centerline{\sl e-mail: \tt 
ney@math.wisc.edu}
\vskip 3mm
\vskip 3mm
\centerline{\bf Abstract}
\bigskip
\bigskip

{\it Consider a finite Abelian group $(G,+)$, with $|G|=p^r$, $p$ a
prime number, and $\varphi: G^\N \to
G^\N$ the cellular automaton given by $(\varphi x)_n=\mu x_n+\nu
x_{n+1}$ for any $n\in \N$, where $\mu$ and $\nu$ are integers
relatively primes to $p$.  We prove that if $\P$ is a translation
invariant probability measure on $G^\Z$ determining a chain with
complete connections and summable decay of correlations, then for any
${\underline w}= (w_i:i<0)$ the Ces\`aro mean distribution
$\displaystyle {\cal M}_{\P_{\underline w}} =\lim_{M\to\infty} {1\over
M} \sum^{M-1}_{m=0}\P_{\underline w}\circ\varphi^{-m}$, where
$\P_{\underline w}$ is the measure induced by $\P$ on $G^\N$
conditioning to $\underline w$, exists and satisfies ${\cal
M}_{\P_{\underline w}}=\lambda^\N$, the uniform product measure on
$G^\N$. The proof uses a regeneration representation of $\P$.}

\medskip
\ni {\bf AMS Classification:} 60K35, 82C, 60K05, 60J05, 58F08.

\ni{\bf Keywords:} group automata, regeneration, chains with complete connections. 

\bigskip
\ni{\bf 1.-- Introduction and main results.}
\numsec=1\numfor=1\bigskip
\bigskip

Let $(G,+)$ be a finite Abelian group with 
$q=p^r$ elements, being $p$ a prime number.  We put 
$\lambda=(q^{-1},..., q^{-1})$ the uniform measure on the
group.  In this paper we study the measure evolution under the
dynamics of 
the cellular automaton $\varphi:G^\N\to G^\N$, given by
$(\varphi x)_n=\mu x_n+ \nu x_{n+1}$ for $n\in\N$, where $\mu$ and
$\nu$ are integers relatively primes to
 $p$ ($\ell g $ means $g+...+g$ $\ell$--times).  
The uniform product measure
$\P=\lambda^\N$ is $\varphi$--invariant,  
$\P\circ\varphi^{-n}=\P$, but any other product measure $\P=\pi^\N$,
with $\pi\neq \lambda$, is not $\varphi$-invariant. Moreover, even in
the simplest case $G=\{0,1\}$ and $+$ the mod 2 sum, the limit of the
marginal distribution, $\lim\limits_{m\to\infty}\P \{ (\varphi^mx)_0=g
\}$ with $g \in G$, does not exist. 
The reason is that for $m=2^k$,
$(\varphi^mx)_0 = x_0+x_m$ (the other terms sum an even number of
times and do not contribute to the sum) has probability
$p^2+(1-p)^2$ to be $0$,
while for $m=2^{k}-1$ this probability converges to $1 \over 2$
because $(\varphi^mx)_0=\displaystyle\sum_{\ell=0}^m x_\ell$.

Alternatively we can study the Ces\`aro mean distribution
$$
{\cal M}_\P\,\doteq\,\lim_{M\to\infty} {1\over M} \sum^{M-1}_{m=0}
\P\circ\varphi^{-m}
$$ 
for a class of initial distributions $\P$ on $G^\N$. In the above
display and in the sequel $\doteq$ means ``it is defined by''.

 Let $-\N^* =\{-i:i\in\N\setminus\{0\}\}$ and $\N^*=\N\setminus\{0\}$.
Let $\P$ be a translation invariant probability measure on $G^\Z$. 
For $\underline w \in G^{-\N^*}$ let $\P_{\underline w}$
be the measure on $G^\N$ induced by the conditional probabilities as 
follows. For
any $m\ge 0$ and $g_0,\dots,g_m\in G$, define
$$
\P_{\underline w}\{x_0=g_0,\dots,x_m=g_m\}\,
\doteq\, \P\{x_0=g_0,\dots,x_m=g_m\,\vert\,
x_{i}=w_{i}, i<0\}.
$$
We
say that $\P$ has {\sl complete connections} if it satisfies
$$
\hbox{ $\forall g_0 \in G, \ \forall \underline w \in G^{-\N^*}$, 
\quad $\P_{\underline w}\{x_0=g_0\}>0. $ }\Eq(50)
$$
For any $m\ge 0$ define
$$ 
\gamma_m\, \doteq\, \sup\left\{\left |{\P_{\underline w}\{ x_0=g
\} \over \P_{\underline v}\{ x_0=g \}} -1 \right | : g \in G, \underline v,
\underline w \in G^{-\N^*}, 
v_{i}=w_{i},i\in[-m,-1]\right\}.
$$ 
We say that $\P$ has {\sl summable decay} if
$$
\sum_{m=0}^{\infty} \gamma_m<\infty. \Eq(60)
$$
This is a uniform continuity condition on $\P_{\underline w}(g)$ as a function
of $\underline w$.

\medskip

The Ces\`aro limits has been already studied for the mod 2 sum
automaton and other classes of permutative cellular automata in [L]
and [MM].  In these papers it is computed mainly for Bernoulli
measures, and in [MM] only the one site Ces\`aro limit is computed for
a Markov measure.  In the mod 2 case the limit is uniformly
distributed, but for some permutative cellular automata the Ces\`aro
mean exists but it is not necessarily uniform. In [FMM] the
Athreya-Ney regeneration times representation of $r$-step Markov chain
was used to show the convergence of the Ces\`aro mean of the group
automata starting with these Markov chains to the uniform Bernoulli
measure.

In this paper we generalize these results for the group automaton
$\varphi$ and initial measures with complete connections and summable
decay.  

\bigskip
\proclaim \Theorem (1).
Let $(G,+)$ be a finite Abelian group with $|G|=p^r$, being $p$ a prime
number.  Let $\P$ be a translation 
invariant probability measure on
$G^\Z$ with complete connections and summable decay.  
Let 
$\varphi:G^\N\to G^\N$ be the cellular automaton, given by
$(\varphi x)_n=\mu x_n+ \nu x_{n+1}$ for $n\in\N$, where $\mu$ and
$\nu$ are integers relatively primes to $p$.  
Then for all
$\underline w \in G^{-\N^*}$ the Ces\`aro mean distribution ${\cal
M}_{\P_{\underline w}}$ exists and verifies ${\cal M}_{\P_{\underline
w}}= \lambda^\N$, the product of uniform measures on $G$. 

\bigskip 

There are two main elements in the proof: regeneration times
and distribution of Pascal triangle coefficients mod $p$.

\bigskip
\noindent{\bf 2.-- Regeneration times for the initial measure.} 
\numsec=2\numfor=1
\bigskip

We show that under the conditions of Theorem \equ(1), for all $\underline w\in
-\N^*$ we can jointly construct a random sequence $\underline x=(x_i:i \in \N)
\in G^\N$ with distribution $\P_{\underline w}$ and a random subsequence
$(T_i:i\in\N^*)\subseteq \N$ such that $(x_{T_i}:i\in \N^*)$ are iid uniformly
distributed in $G$ and independent of $(x_i:i\in \N \setminus
\{T_1,T_2,\dots\})$; furthermore $(T_i:i\in\N^*)$ is a stationary renewal
process with finite mean inter-renewal time independent of $\underline w$. A
consequence of the construction is that the random vectors (of random lenghts)
$((x_{T_i},\dots, x_{T_{i+1}-1}): i\ge 1)$ are iid. 

Our regeneration approach shares results with Berbee (1987) and Ney and
Nummelin (1993). The construction is simple: the probability space is
generated by product of iid uniform (in $[0,1]$) random variables. It works as
the well known construction and simulation of Markov chains as a function
of a sequence of uniform random variables (see, for instance Ferrari and
Galves (1997)). Bressaud, Fern\'andez and Galves (1998) construct a coupling
using these ideas to show decay of correlations for measures with infinite
memory.

For $\underline w \in G^{-\N^*}$ and $g \in G$ denote
$$
P(g\vert \underline w)\doteq \P\{x_0=g\vert x_i=w_i,
i\le -1\}.
$$

Let  
$$
a_{-1}(g\vert \underline w)\,\doteq \,
\inf \{ P(z\vert \underline v)\ : \  {\underline v} \in G^{-\N^*},\, z \in 
G \}. \Eq(g1)
$$ Actually $a_{-1}$ depends neither on $g$ nor on $\underline w$; we keep the
dependence in the notation for future (notational) convenience.  Since the
space $G^{-\N^*}$ is compact and
$\P$ has summable decay, the infimun in \equ(g1) must be attained by a
$g^0\in G$ and a $\underline w^0\in G^{-\N^*}$. Hence,
$$
a_{-1}(g\vert \underline w)\, =\,P(g^0\vert \underline w^0) \,>\, 0,
$$
because $\P$ has complete connections.
For each $k \in\N$, $g \in G$ and ${\underline w} \in G^{-\N^*}$
define 
$$
a_k(g \vert {\underline w}) \doteq 
\inf\{P(g\vert w_{-1},\dots, w_{-k}, {\underline z}) \ : \ 
{\underline z} \in G^{-\N^*} \},
$$
where $(w_{-1},\dots, w_{-k}, {\underline z}) = (w_{-1},\dots, w_{-k},
z_{-1}, z_{-2},\dots)$. Notice that $a_0(g\vert \underline w)$ does not depend
on $\underline w$. 
Let 
$$
b_{-1}(g\vert \underline w) 
\,\doteq\, a_{-1}(g\vert \underline w) 
,
$$ 
for $g\in G$. For $k\ge 0$,
$$ 
b_k(g \vert {\underline w}) \,\doteq\, a_k(g \vert {\underline
w})-a_{k-1}(g \vert {\underline w}).
$$

We construct disjoint intervals $B_k(g \vert {\underline w})$ for
$g\in G$, $k\ge -1$, contained in $[0,1]$, of Lebesgue measure $b_k(g
\vert {\underline w})$ respectively, disposed in increasing order with
respect to $g$ and $k$: $ B_{-1}(0|\underline w), \dots,B_{-1}(q-1|\underline w),
B_0(0|\underline w),\dots,B_0(q-1|\underline w),B_1(0|\underline
w),\dots, B_1(q-1|\underline w), \dots$, with no intersections (we
have enumerated $G$ by $\{0,...,q-1\}$). The construction guarantees
$$
\Bigl|\bigcup_{k\ge -1} B_k(g \vert {\underline
w})\Bigr|\,=\,P(g|\underline w)
$$
and
$$ 
\Bigl|\bigcup_{g\in G}\,\bigcup_{k\ge -1} B_k(g
\vert {\underline w})\Bigr|\,=\, 1.
$$ 
(All the unions above are disjoint.)

Let $\underline U= (U_n: n \in \Z)$ be a double infinite sequence of
iid random variables uniformly distributed in $[0,1]$. Let
$(\Omega,{\cal F},\P)$ be the probability space induced by these
random variables. For each $\underline w\in G^{-\N^*}$ we construct
the random sequence $\underline x$ with distribution $\P_{\underline
w}$ in $\Omega$, as a function of $\underline U$, recursively: for
$n\in \N$
$$ 
x_{n} \,\doteq \,\sum_{g\in G} g \Bigl[ \sum_{\ell\ge -1} \one\{ U_n \in
B_\ell(g \vert x_{n-1},\dots,x_0,{\underline w}) \}\Bigr].
$$ 

For $\ell\ge -1$ let
$$
B_\ell(\underline w) \doteq 
\bigcup_{g\in G}  B_\ell(g|\underline w).
$$
Notice that neither $B_{-1}(g|\underline w)$ nor $B_{-1}(\underline
w)$ depend on $\underline w$. Furthermore
$$
{|B_{-1}(g|\underline w)|\over|B_{-1}(\underline w)|} \, =\, |G|^{-1}.\Eq(gm1)
$$
For $k \in \N$ let
$$
a_k \doteq \min_{\underline w} \ \left \{ \displaystyle  \sum_{g\in G} 
 a_k(g \vert \underline w)\right\}.
$$
This is a non-decreasing sequence and satisfies
$$
[0,a_k]\;\subset \;\bigcup_{\ell=-1}^k 
B_\ell(\underline w), \Eq(oak)
$$ 
independently of $\underline w \in G^{-\N^*}$.
\bigskip
\proclaim\Lemma(71). In the event $\{U_n\le a_k\}$ for $n\in \N$ 
we only need to look at
$x_{n-1},\dots,x_{n-k}$ to decide the value of $x_n$. More precisely,
for $\underline v\in G^\Z$ such that $v_i=w_i$ for $i\le -1$,
$$
\eqalign{
\P_{\underline w}\{x_n=g \,\vert\, &U_n\le a_k,\,x_{n-1}=v_{n-1}, \dots, 
x_0=v_0\}\,\cr
&= \,\P_{\underline w}\{x_n=g
\,\vert\, U_n\le a_k,\,x_{n-1}=v_{n-1}, \dots, x_{n-k}= v_{n-k}\}.}
$$

\noindent {\bf Proof.} Follows from \equ(oak). \square
\medskip

Define times
$$
\eqalign{
T_1&\doteq \min \{n\ge 0: \,U_{n+j} \le a_{j-1}, \
j\ge 0 \},\cr
T_i &\doteq \min \{n>T_{i-1}: \,U_{n+j} \le a_{j-1}, \
j\ge 0 \}, }
$$ 
for $i > 1$.  
\bigskip

Let $\bn$ be the counting measure on $\N$ induced by $(T_i:i\ge 1)$:
for $A\subset \N$ and $n\in\N$, 
$$ \bn(A) \,\doteq\, \sum_{i\ge 1} \one\{T_i\in A\},\ \ \
\bn(n)\,\doteq\,\bn(\{n\}).
$$
Notice that the definitions of $(T_i:i\ge 1)$ and $\bn$ depend only on $(U_n:n\in \Z)$
and do not depend on~$\underline w$.

\proclaim\Lemma(70). The distribution of the counting measure $\bn$
corresponds to a stationary renewal process.

\noindent{\bf Proof. } We will construct a stationary renewal process
$\bm$ in $\Z$ whose projection on $\N$ is $\bn$.
For $k\in\Z$, $k'\in\Z\cup\{\infty\}$, define
$$
H[k,k'] \,\doteq\,\cases{  \{U_{k+\ell}\le a_{\ell-1},
\ell=0,\dots,k'-k\},&if $k\le k'$\cr
&\cr
\hbox{``full event''},&if $k>k'$ }
$$
With this notation,
$$
\bn(n) = \one \{H[n,\infty]\},\ \ n\in\N. \Eq(500)
$$
We construct a double infinity counting process $\bm$ using the
variables $(U_n:n\in \Z)$ by 
$$
\bm(n) \,\doteq\, \one \{H[n,\infty]\},\ \ n\in\Z.
$$ 
By construction, the distribution of $\bm$ is translation
invariant, hence $\bm$ is stationary. Furthermore, by \equ(500) it
coincides with $\bn$ in $\N$: $\bm(K)=\bn(K)$ for $K\subset \N$. 
Define $T_i$ for $i\le 0$ as the ordered time-events of $\bm$ in the
negative axis. 

The (marginal) probability of a counting event at time $n\in\Z$, is
given by
$$ 
\P\{\bm(n)=1\}
= \P\{U_{n+j} 
\le a_{j-1}, \ j\ge 0\}
\,=\, a_{-1}\, a_0\, a_1\,\cdots 
\,\doteq\, \beta,
$$
and it is independent of $n$. 
We first show that under the hypothesis of summability of $\gamma_k$,
$\beta $ is strictly positive. 
For any $g \in G$, $w_{-1},\dots, w_{-k} \in G$ 
and ${\underline z}, {\underline v}
\in  G^{-\N^*}$
$$
\left | {\P\{g\vert  w_{-1}\dots w_{-k}, {\underline z}\} \over 
\P\{g\vert  w_{-1}\dots w_{-k}, {\underline v}\}}
-1 \right | \le \gamma_k,$$
therefore
$$
\inf \ \{ \P\{g \vert  w_{-1}\dots w_{-k}{\underline z}\} \ : \
{\underline z} \in G^{-\N^*} \}\ge 
(1-\gamma_k) \P\{g\vert  w_{-1}\dots w_{-k} {\underline v}\}.
$$
Summing over $g\in G$ and taking minimum on the set
$\{w_{-1},\dots,w_{-k}\}$ we conclude that
$$a_k \ge 1-\gamma_k.$$ 
Since $\displaystyle 
\sum_{k\ge 0} \gamma_k < \infty$ we deduce that 
$\displaystyle 
\sum_{k\ge 0} (1-a_k) < \infty$ and henceforth 
$\beta >0$.

\medskip

We show now that $\bm$ is a renewal process on $\Z$.
The event $\{\bm(n) = 1\}$ depends only on $(U_k:k\ge n)$, that is,
$(T_i:i\in\Z)$ are stopping times for the process $(U_{-k}:k\in \Z)$.
Since for $k<k'<k''\le \infty$, 
$$
H[k,k'']\cap H[k',k'']\,\, =\,\, H[k,k'-1]\cap H[k',k''],
$$ we have that for any finite set $A=\{k_1,\dots,k_n\}$ with
$k_1<\dots<k_n<k'$ and for any sequence $(m_\ell:\ell>k')$ with
$m_\ell\in\{0,1\}$,
\smallskip
$$
\eqalign{
\P\big\{&\bm(A)=n\,\big\vert\, \bm(k')=1, \bm(\ell)= m_\ell,
\,\ell>k'\big\}\cr
&=\,\P\big\{\displaystyle \bigcap_{i=1}^nH[k_i,k'-1]\,\big\vert\, \bm(k')=1\big\}
\phantom{\sum_o^p}\cr 
&= \,\prod_{i=1}^n\P\{H[k_i,k_{i+1}-1]\},\cr
}\Eq(epb)
$$ 
where $k_{n+1}\doteq k'$. The computation above could be done
because $\P\{\bm(k')=1\}=\beta>0$.  Display \equ(epb) means that given a
counting event at time $k'$, the distribution of the counting events
for times less than $k'$ does not depend on the events after
$k'$. This characterizes $\bm$ as a renewal process.  Since the
density $\beta$ is positive, $T_1$, the residual time is a honest
random variable, and for $i\neq1$,
$\E(T_{i+1}-T_i)=\beta^{-1}<\infty$.\ \ \square

\bigskip
\proclaim\Lemma(72). The variables $(x_{T_i}:\,i\ge 0)$ are iid
uniformly distributed in $G$.

\noindent {\bf Proof.} Let us show that the marginal distribution of
$x_{T_i}$ is uniform in $G$. Since times $(T_i:i\in \N^*)$ 
are finite almost surely:
$$
\eqalign{
\P\{x_{T_i}=g \}
&= \sum_{n\in \N}\P\Bigl \{U_n \in \displaystyle \bigcup_{\ell \ge -1}
B_\ell(g\vert \underline w) \,,\,T_i = n\Bigr \}\cr
&= \sum_{n\in \N}\P \{U_n \in B_{-1}(g\vert\underline w)\,\vert\, U_n\in
B_{-1}(\underline w)\}\, \P\{T_i = n\} \cr
&= |G|^{-1}.}  
$$ 
The second identity follows because $\{T_i=n\}$ is the intersection
of $\{U_n \in  B_{-1}(\underline w)\}$ with events depending on variables
$(U_{n+\ell}, \, \ell\neq 0)$ which are independent of $U_n$. The
third identity follows from \equ(gm1). The
same computation shows that for any $K\subset \N$, $(i(k):k\in
K)\subseteq \N$, and $(g_k:k\in K)\subseteq G^{K}$
$$
\P\{x_{T_{i(k)}} = g_k,\, k\in K\} = |G|^{-|K|},
$$ 
so that $(x_{T_{i(k)}}:k\in K)$ are iid in $G$.  The reason why
the above computation works is that in the event $\{T_i=n\}$,
$U_{n+1}\le a_0$, hence $x_{n+1}$ does not depend on the past. Since for
all $j\ge 1$, $U_{n+j}\le a_{n+j-1}$, $x_{n+j+1}$ only depends on
$x_{n+1},...,x_{n+j}$.
\square

\medskip\medskip
\ni{\bf 3.-- A renewal Lemma.}
\medskip\medskip
\numsec=3 \numfor=1

In this section we show that a stationary discrete-time renewal
process on $\N$ has high probability to visit sets with many points.

\proclaim \Lemma (4).  Let $\bn$ be a stationary renewal process with
finite inter-renewal mean. Then for all $A\subset \N$,
$$
\P\{\bn(A)=0\} \le \vep(|A|)
$$ 
with $\vep(n)\to 0$ as $n\to\infty$. Also, $\vep :\N \to \R$ can be
chosen to be decreasing.


\noindent{\bf Proof.} We are going to prove that for all $\vep>0$
there exists $n_0$ such that for any finite set $A\subset \N$ with
$|A|>n_0$,
$$
\P\{\bn(A)=0\} \le \vep. \Eq(221) 
$$ 
We start with some known facts of renewal theory.  Let $T_i$ be the
renewal times and $\beta = 1/\E(T_{i+1}-T_i)$ for some $i\ge 1$ (and
hence for all $i\ge 1$). Since the inter-renewal distribution has a
first moment finite, the key renewal theorem holds: we have
$$
\lim_{n\to\infty}\P\{\bn(n) =1\,\vert\, \bn(0)=1\}\,=\,  \beta. \Eq(220)
$$
Let $S_n\doteq T_{\bbn(n)+1}-n$ be the residual time (over jump) at $n$,
where we have denoted by $\bbn(n)=\bn([0,n])$,
and let for $k\ge 0$
$$
\beta_k=\P\{\bn(k)=1\,|\, T_1=0\},\ \ \ F(k) = \P\{T_2-T_1>k\},\ \ \ 
F_n(k) = \P\{ S_n >k \, | \, T_1=0 \}. \Eq(221)
$$
Now we have
$$
F_n(k)\,=\, \sum_{j=0}^n F(j+k) \beta_{n-j}\,\le\, \overline F(k), \Eq(222)
$$
where
$$
\overline F(k)\,\doteq\, \sum_{j=k}^\infty F(j)\,\to \,0
$$
as $k\to\infty$ because we are assuming that the inter-renewal time
has a finite mean.

For any subset $B\subset A$ we have
$$
\P\{\bn(A)=0\}\,\le\,\P\{\bn(B)=0\}.
$$
For any $A$ with $|A|=n$ and any $1<\ell<n$, there exists a set
$$
\{b^n_1,\dots,b^n_\ell\} \doteq B^n_\ell\subset A
$$
with
$$
\Bigl[{n\over \ell}\Bigr]\,\le\, b^n_{j+1}-b^n_j, \ \ \ j=1,\dots,\ell-1,\Eq(224)
$$
where $[x]$ is the largest integer in $x$. The choice of 
$\{b^n_1,...,b^n_\ell\}$
depends on $A$ but $\ell$, and \equ(224) hold uniformly for all $A$
with $|A|=n$.

Let $\vep>0$ and take any $0<\delta<\beta$. Take $n_0$ such that
$\beta_n>\delta$ for $n>n_0$. Let  $n>\ell n_0$ and define 
$$ 
\Gamma^n_j \,\doteq\, \{S_{b^n_j}\,\le [n/\ell] - n_0\}
$$
the event ``the over jump of $b^n_j$ does not superate $[n/\ell]-n_0$''.
Let
$$
\Theta^n_j\,\doteq\,
\bigl\{\bn(b^n_j-b^n_{j-1}-S_{b^n_{j-1}})\,=\,0\bigr\},
\Lambda^n_j\,\doteq\,
\bigl\{\bn(b^n_j)=0\bigr\}
$$
the events ``starting at the over jump of $b^n_{j-1}$, $b_j^n$ is not
hit'' and ``$b_j^n$ is not hit'' respectively.


\ni From \equ(222) we get for $2\le j \le \ell$
$$ \P \{\Gamma^n_j \,|\,\Gamma^n_{j-1}\} \ge 
\big(1- \overline F([n/\ell]-n_0)\big). \Eq(226)$$

\ni Then
$$
\eqalign{
\P\{\bn(A) =0\} \,&\le\, \P\{\bn(B^n_\ell)=0\}
= \P\{\Lambda_1^n\cap .... \cap \Lambda_{\ell}^n \} \cr
&\le \P\{\Lambda_1^n\cap \Gamma^n_1\} 
\P\{\Lambda_2^n\cap .... \cap \Lambda_{\ell}^n \,|\,\Gamma^n_{1} \} +
1-\P\{\Gamma_1^n\} \cr
&\le \prod_{j=1}^{\ell} \P\{\Theta_{j}^n \,|\,\Gamma^n_{j}\}
+\sum_{j=1}^{\ell-1} (1 - \P\{\Gamma^n_{j} \,|\, \Gamma^n_{j-1}\} ) \cr
&\le (1-\delta)^{\ell} + (\ell -1) \overline F([n/\ell]-n_0) +
\P\{T_1 > [n/\ell]-n_0 \} \cr
}\Eq(225)
$$
since $\beta_n > \delta$ for $n > n_0$ and \equ(226).
Now choose $\ell$ so that $(1-\delta)^\ell<\vep/3$, then $n$ so that
$\overline F([n/\ell]-n_0)<\vep/3(\ell-1)$ and 
$\P\{T_1 > [n/\ell]-n_0\} \le \vep / 3$, to conclude
$$
\P\{\bn(B^n_\ell)=0\}\,\le\,\vep
$$
for sufficiently large $n+\ell$. \square

\bigskip
\bigskip
\ni{\bf 4.-- Convergence of Ces\`aro limit.}
\bigskip
\bigskip
\numsec=4 \numfor=1

For proving this theorem we shall need some results concerning 
walks of variables
determining a chain with complete connections. 
In this purpose let 
us introduce some notation.
First $R=(r_k:k\in\N)$ denotes an increasing sequence in $\N$. 
We put  
$R_n=(r_k:k\le n)$. For any subsequence 
$\overline R=(\overline r_k: k \in \N)$ of $R$ we define the index function by
$f_{\overline R}(k)=\ell $ if $\overline r_k =r_\ell$. We also set 
$n(\overline R)=|\overline R\cap R_n|$.
Let $a^R=(a^R_r:r\in R)$ be a sequence of non-negative integers. They
define maps  $\psi^R_r:G\to G$ such that $\psi_r^R(g)=a_r^R \ g=g+...+g  $ \
$a_r^R$ times, for any $r \in R$.
We associate to it the following 
sequence of random variables taking values in $G$,
$$
S^R_n =\sum\limits_{r\in R_n} a_r^R \ x_r,\ n\in\N. 
$$
We will distinguish  the following subsequence 
$$R^* \doteq R^*(a ^R) = \{ r\in R \ : \ a_r^R\neq 0 \ mod \ p  \}.$$

\medskip
\ni {\bf Remark.} Since $(G,+)$ is a finite Abelian group with 
$|G|=p^r$, $p$ a prime number, then the function $\psi(g)=a\ g $, 
where  $a\in \N$, is one-to-one whenever $a\neq 0 \ mod \ p$.
\medskip

Let $J\subseteq \N$ be a finite set. Consider 
a finite family of sequences $R^J=(R^j:j\in J)$. Associated to each
sequence there is a sequence of non-negative 
integers $a^{R^j}=(a_r^{R^j}:r \in R^j)$ and the corresponding
set of mappings
$\psi^{R^j}= (\psi^{R^j}_r:r\in R^j)$.
As before we consider the sequences
$R^{j*}\doteq R^*(a^{R^j})$ for 
$j\in J$. Let $\tilde R^J=(\tilde R^j:j\in J)$ 
be a family of subsequences verifying the following conditions:
\medskip
\item{(H1)} $\tilde R^j\subseteq R^{j*}$ for any $j \in J$, 
\item{(H2)} $\tilde R^j \cap \tilde R^i=\emptyset$ 
if $i \neq j$ in $J$, 
\item{(H3)} if $r \in \tilde R^j\cap R^k $ for $k<j$ in $J$, then
$a_r^{R^k}=0 \ mod \ p$. 
\medskip
\ni We set 
$$\tilde n(\tilde R^J)=\min\{ n(\tilde R^j):j\in J\}$$
and
$$\tilde n(R^J)=\max\{\tilde n (\tilde R^J):\tilde R^J 
\hbox{ verifying } (H1), (H2), (H3) \}.$$
\medskip
The proof of  Theorem \equ(1) is based upon the   following result.

\bigskip
\proclaim \Lemma (401). 
Let $\P$ be a translation invariant measure on $G^\Z$ 
with complete connections such that
$\displaystyle\sum_{m\ge 0} \gamma_m < \infty$, and let $\underline w \in G^{-\N^*}$.
Then
\item{(a)} $\exists \ \vep_1:\N \to \R$, a decreasing function with 
$\vep_1(n)\to 0$ if $n \to \infty$,  
such that  for any increasing sequence $R$ in $\N$ and any
sequence of non-negative integers $a^R$ it is verified 
$$
\left |\P_{\underline w}\{ S^R_n=g\}-q^{-1} \right | \le \vep_1(n(R^*)), \ \hbox{ for any } 
 n\in \N, g\in G.
$$
\item{(b)} Let $J \subset \N$ be  finite.  Then there is a decreasing function
$\vep_J:\N \to \R$ with  $\vep_J(n)\to 0$ if $n\to \infty$, 
such that  
for any set of sequences $R^J=(R^j:j\in J)$ and any  family of
non-negative integers   
$(a^{R^j}:j\in J)$, it is verified 
$$\left |\P_{\underline w}\{S^{R^j}_n=g_j, \hbox{ for }  j\in J\} -q^{-|J|}\right |\le 
\vep_J(\tilde n(R^J)), \ 
\hbox{ for any }n\in \N, (g_j:j\in J)\in G^{J} .
$$
\square

\bigskip
Before begin the proof of Lemma \equ(401) we include a useful arithmetic property. We include
a proof for completeness.
For $(G,+)$ a finite  Abelian group with $|G|=p^r$, where
$p$ is a prime number, consider the following system 
of equations (S):
\bigskip
$$\matrix{
(1)&a_{11}g_1 &+&a_{12}g_2&+&...&+&a_{1\ell}g_{\ell}=0 \cr
(2)&a_{21}g_1 &+&a_{22}g_2&+&...&+&a_{2\ell}g_{\ell}=0 \cr
   &          & &         & &\vdots   & &                    \cr
(\ell)&a_{\ell1}g_1 &+&a_{\ell 2}g_2&+&...&+&a_{\ell \ell}g_{\ell}=0 \cr
}
$$
such that 
$$
(H') \qquad a_{ij} \in \N, \ a_{ii}\neq 0 \ mod \ p, \ a_{ij}=0 \
mod \ p  \hbox{ if } i < j.
$$ 
Denote $a_{ii}=k_ip+s_i$ with $s_i \in \{1,...,p-1\}$
and $a_{ij}=c_{ij}p$ for $i<j$.
\bigskip

\proclaim \Lemma (402).  
The system (S) has unique solution $g_1=g_2=...=g_{\ell}=0$.

\smallskip
\ni {\bf Proof.}
First of
all we will prove that if $g_1,...,g_{\ell}$ are solutions of
(S) and for some $1 < s \le r$,  $p^sg_i=0$, $i \in \{1,...,\ell\}$,
then $p^{s-1}g_i=0$ for $i \in \{1,...,\ell\}$. 
We prove this property by induction on $\{1,...,\ell\}$.
First consider equation (1),
$$(k_1p+s_1) g_1+\sum_{j=2}^{\ell}c_{1j}pg_j=0.$$
If we add the equation $p^{s-1}$ times we obtain,
$$k_1p^sg_1+s_1p^{s-1}g_1+\sum_{j=2}^{\ell}c_{1j}p^sg_j=0,$$
then $s_1p^{s-1}g_1=0$. Since the product by $s_1$ defines 
a 1-to-1 map we conclude that $p^{s-1}g_1=0$.
Let us  continue with the induction assuming that $p^{s-1}g_1=0$,
$p^{s-1}g_2=0,...,p^{s-1}g_t=0$, for $1\le t < \ell$, and we prove that
$p^{s-1}g_{t+1}=0$.
\medskip
Adding $p^{s-1}$ times equation $t+1$ we get
$$\sum_{j=1}^ta_{t+1,j}p^{s-1}g_j+(k_{t+1}p+s_{t+1})p^{s-1}g_{t+1}+
\sum_{j=t+2}^{\ell}c_{t+1,j}p^{s}g_j=0.$$
Therefore, using the induction hypothesis we obtain
$s_{t+1}(p^{s-1}g_{t+1})=0$ and henceforth $p^{s-1}g_{t+1}=0$.
\medskip
To conclude we use last property recursively
beginning from the fact that $p^rg_i=0$ for
any $i \in \{1,...,\ell\}$. \square

Hence the transformation $A:G^\ell \to G^\ell,$ $A \vec g=\vec h$, with $\vec g,
\vec h \in G^\ell$ and matrix $A$ verifying condition (H') is a one-to-one
and onto transformation. In what follows we identify $\P_{\underline w}$
with $\P$.
\bigskip
\ni {\bf Proof of Lemma \equ(401).}
\bigskip

\noindent a) For any increasing sequence $R= (r_k:k\in\N)$ we put
$$
\tau^{R} = \inf\{k\in \N: \bn(r_k)=1\}, 
\hbox{ where } \infty=\inf\phi, 
$$
the first time that some element of the sequence $R$ belongs to the
renewal process $\bn$ introduced in Section 2. Consider $R^*$ the
subsequence corresponding to  mappings $\psi_r^R$ 
such that
$a_r^R \neq 0 mod p$. We
denote $n^*=n(R^*)$, $\tau^*=\tau^{R^*}$ and $f=f_{R^*}$ the
corresponding index function.  First we prove
$$
\P\{S_n^R=g\vert \tau^{*}\le n^*\}=q^{-1} .
$$
To see that write
$$
\eqalign{
&\P \left \{ S_n^R=g, \tau^{*}\le n^* \right \}= \sum_{k=o}^{n^*} 
\P\left \{S_n^R=g, \tau^{*}=k \right \}\cr
&=\sum_{k=0}^{n^*} \P \left\{\sum_{i=0}^{f(k)-1}\psi_{r_i}(x_{r_i})+
\psi_{r_{f(k)}}
(U_{r_{f(k)}})+
\sum_{i=f(k)+1}^{n}\psi_{r_i}(x_{r_i})=g, \tau^{*}=k\right\} \cr
&=\sum_{k=0}^{n^*}\sum_{g_1, g_2 \in G} 
\P \Bigl \{  \sum_{i=0}^{f(k)-1}\psi_{r_i}(x_{r_i})=g_1, 
U_{r_{f(k)}}= \psi^{-1}_{r_{f(k)}}(g-g_1-g_2),  \cr
&  \hskip 4cm \qquad\qquad\qquad\qquad  
\sum_{i=f(k)+1}^{n}\psi_{r_i}(x_{r_i})=g_2, \tau^{*}=k \Bigr \}  \cr}$$
$$\eqalign{
&=q^{-1}\sum_{k=0}^{n^*}\sum_{g_1, g_2 \in G} 
\P\left\{\sum_{i=0}^{f(k)-1}\psi_{r_i}(x_{r_i})=g_1, 
\sum_{i=f(k)+1}^{n}\psi_{r_i}(x_{r_i})=g_2, \tau^{*}=k\right\}  \cr
&=q^{-1}\sum_{k=1}^{n^*} 
\P \left \{ \tau^{*}=k \right \}= 
q^{-1} \P \{\tau^{*}\le n^* \}. \cr
}
$$
Where in the last equalities we have used that $U_{r_{f(k)}}$ is independent
of variables $(x_n:n\neq  r_{f(k)})$ when $\tau^*=k$.
Then,
$$
\P\{S^R_n=g\} =q^{-1}\P\{\tau^{*}\le n^*\} +
\P\{ S^R_n =g, \tau^{*} >n^*\}
$$
and
$$
 \P\{S^R_n=g\} - q^{-1}  = - q ^{-1} \P\{\tau^{*} > n^*\} +
\P\{ S^R_n =g, \tau^{*} >n^*\}.
$$
Using Lemma \equ(4)we get
$$
\left |\P\{S^R_n=g\} - q ^{-1} \right |\le 2 \P \{\tau^{*}>n^*\}\le 
2 \vep(n^*+1). 
$$
\ni b) 
Let $R^J=(R^j:j\in J)$ be a family of sequences,
$(a^{R^j}:j \in J)$ be the family of non-negative sequences, 
$(\psi^{R^j}:j\in J)$ be the corresponding family of 
mappings and $\tilde R^J$ be a family of subsequences verifying 
conditions (H1), (H2), (H3).
Denote by $f_j=f_{\tilde R^j}$ and $\tau_j=\tau^{\tilde R^j}$ for any $j\in J$.
Fix $n \in \N$ and 
put $\tilde n=\tilde n(\tilde R^J)$.  

Take a vector
$\vec k=(k_j:j \in J) \in \{1,...,\tilde n\}^J$.
On the set $\{\tau_j=k_j:j\in J\}$ we define the random
variables 
$$\rho_j(\vec k,n,\underline U)=
\sum_{i \in J} \ \one\{\tilde r_{k_i}^i\in R^j_n\} \ \psi_{f_i(k_i)}(U_{f_i(k_i)}),
\hbox{ for } j \in J.$$
Consider $(g_j':j\in J)\in G^J$. From hypothesis
(H1), (H2), (H3) the system of linear equations
$\rho_j(\vec k,n,\underline U)=g'_j$, $j\in J$, defines  
a system of type (S). Then, by Lemma \equ(402),
there is a unique  $(g_j'':j\in J) \in G^J$ such that
$$\rho_j(\vec k,n,\underline U)=g'_j, \ j\in J\ 
\Leftrightarrow \ U_{f_j(k_j)}=g''_j, \ j \in J. \Eq(403)$$
Let 
$T(\vec k)=(\displaystyle\bigcup_{j\in J} R_n^j)\setminus \{f_j(k_j):j\in J\}.$
It is easy to see that variables $(S_n^{R^j}: j\in J)$ on 
$\{\tau_j=k_j:j\in J\}$ can be written as
$$S_n^{R^j}=\sum_{r \in T(\vec k)\cap R^j_n} \ \psi_r(x_r) \ + \
\rho_j(\vec k,n,\underline U).$$ 
Therefore,
$$\eqalign{
& \P\{S^{R^j}_n=g_j, \tau_j =k_j, \hbox{ for } j\in J\}=\cr 
&\sum_{h_r\in G: \ r\in T(\vec k)} 
\P\{\rho_j(\vec k,n,\underline U)  =g_j- \hskip -0.5 cm \sum_{r\in T(\vec k)\cap R_n^j} 
\psi_r(h_r), 
x_r=h_r, 
\tau_j=k_j, \hbox{ for } j\in J,  
r \in T(\vec k)\} \cr
&=\sum_{h_r \in G:\ r \in T(\vec k)}
\P\{ U_{f_j(k_j)}=g''_j, x_r=h_r, \tau_j=k_j, 
\hbox{ for } \ j \in J,
r \in T(\vec k)\},  \cr
}
$$
where $(g_j'':j \in J) \in G^J$ is given by property \equ(403). 
By independence we conclude that
$$\P\{S_n^{R^j}=g_j, \tau_j=k_j, \hbox{ for } j\in J\}=
q^{-|J|} \P\{\tau_j=k_j, \ j \in J\}.
$$
Hence
$$
\P\{S^{R^j}_n = g_j,\hbox{ for }
j\in J, \ \max\limits_{j\in J} \tau_j\le \tilde n\}=
q^{-|J|}\P\{\max\limits_{j\in J} \tau_j\le \tilde n\},
$$
which together with Lemma  \equ(4) allow us   to deduce that
$$\left |\P\{S^{R^j}_n=g_j, \hbox{ for } j\in J\}-q^{-|J|}\right | \le 2 
\P\{\max\limits_{j\in J}\tau_j>\tilde n\}\le 2|J| \vep(\tilde n +1).$$
\square
\bigskip
Now we can give the proof of the main theorem.
\bigskip
\noindent{\bf Proof of Theorem \equ(1).} 
\bigskip
\bigskip
First, let us introduce some notation. The $p$-expansion of $m\in\N$
is $m=\sum\limits_{i\ge 0} m_ip^i$ with $m_i\in\Z_p$. We denote by
${\cal I}(m)=\{i\in\N :m_i\neq 0\}$ its support
and we denote its elements in decreasing 
order, ${\cal I}(m)=\{\delta_{1,m}>...>\delta_{s_m,m}\}$, where
$s_m=|{\cal I}(m)|$. Now put $m^{(i)}=m_{\delta_{i,m}}$, so $m=\sum\limits^
{s_m}_{i=1} m^{(i)} p^{\delta_{i,m}}$.  Observe that $\delta_{1,m}=$
integer part $(\log m)$, where we take $\log m$ in base $p$.

Since $p$ is a prime number the Lucas' theorem [Lu] asserts that
$$\left[ {m\choose k}\right]_p =\left[ \prod\limits_{i\ge 0}{m_i\choose k_i}
\right]_p, $$
where ${r\choose s}=0$ if $r<s$. In particular
$[{m\choose k}]_p>0$ if and only if $k_i\le m_i$ for all
$i\ge 0$. 

Let us return to the automaton $\varphi$. Since $G$ is Abelian, a simple recurrence
implies
$$
(\varphi^mx)_i=\sum\limits_{k\le m} {m\choose k} \mu^{m-k} \nu ^k x_{k+i}.$$
Observe that
this expression has the form of variables
$S_n^R$ defined before. In this case the mapping
has the shape $ {m \choose k } \mu^{m-k}\nu^k \ g$
which is one-to-one if $\left[ {m \choose k } \right ]_p \neq 0$
since $\mu$ and $\nu$ are relatively primes to $p$.
Then our computations are devoted to show that 
we have enough one-to-one mappings.

In order to make the proof more clear 
we shall first prove that the Ces\`aro mean of the marginal distribution
exists and it is uniform, that means  
$$
\pi(g) \doteq\lim\limits_{M\to\infty}{1\over M}\sum\limits^{M-1}_{m=0}
\P_{\underline w}\{(\varphi^mx)_0=g\}\hbox{ exists and verifies }
\pi(g)=q^{-1},\hbox{ for any }g\in G.
$$

Let us fix 
$\alpha\in (0,{1\over 2})$. For
$M>0$ consider the set
$
{\cal R}_M=\{m\le M: |{\cal I}(m)| \ge \alpha\log\log M\}.$ We will
prove that $({\cal R}_M:M\in \N)$ is a sequence of sets of density one,
which means 
${|\{m\le M\}\setminus {\cal R}_M| /  M}
\sobre{\hbox{\rightarrowfill}}{M\to\infty}0$. 
In that purpose we make the 
decomposition \hfill\break
$
\{m\le M\}=\bigcup\limits_{1\le s\le s_M+1} A_{s,M}
$
with 
$$
\eqalign{
&A_{1,M}= \{m\le M: \delta_{1,m}<\delta_{1,M}\},\cr
&A_{s,M}= \{m\le M: \delta_{r,m}=\delta_{r,M}
\hbox{ for } r<s \hbox{ and } \delta_{s,m}<\delta_{s,M}\} \hbox{ for } 
1\le s\le s_M,\cr
&A_{s_M+1,M}= \{M\}. }
$$
Observe that  $|A_{s,M}|=M^{(s)} p^{\delta_{s,M}}$ for $1\le s\le s_M$. 
Take $s^*_M=\sup\{ s:\delta_{s,M}\ge \log\log M\}$. Since
$\delta_{1,M}=$ integer part $(\log M)$, we have $s^*_M\ge 1$. Now,
$$
|\{m\in A_{s,M}: |{\cal I}(m)|\le \alpha\delta_{s,M}\}|\le \sum\limits_{t\le
\alpha\delta_{s,M}} (p-1)^t{\delta_{s,M}\choose t}
$$
$$\le(p-1)^{\alpha\delta_{s,M}}2^{\delta_{s,M}}e^{-2(\alpha-{1\over 2})^2
\delta_{s,M}}.$$
Hence,
$$
|\{m\le M\} \setminus {\cal R}_M|\le\sum\limits_{1\le s\le s^*_M}
(2(p-1)^\alpha)^{\delta_{s,M}}e^{-2(\alpha-{1\over 2})^2
\delta_{s,M}}
+\sum\limits_{s_M^*<s\le s_M} M^{(s)} p^{\delta_{s,M}}+1.
$$
We have
$$
\sum\limits_{s_M^*<s\le s_M} M^{(s)}  p^{\delta_{s,M}}+1\le 
(\log M)^2 +1 .
$$
Take  $\alpha <{p\over 2}(\log(p-1))^{-1}$, then  $p'\doteq 2(p-1)^
\alpha e^{-2(\alpha-{1\over 2})^2} < p$.  Therefore 
$${1 \over M} \sum\limits_{1\le s\le s^*_M}(2(p-1)^\alpha)^{\delta_{s,M}}e^
{-2(\alpha-{1\over 2})^2 \delta_{s,M}} \le {1\over M}
\sum\limits_{1\le s\le s^*_M} p'^{\delta_{s,M}}$$
$$
\le \sum\limits_{1\le s\le s^*_M}\( {p'\over p}\)^{\delta_{s,M}} \le
{p\over p-p'}\({p'\over p}\)^{\log \log M}.
$$
Hence
${|\{m\le M\}\setminus {\cal R}_M| / M}
\sobre{\hbox{\rightarrowfill}}{M\to\infty}0$. 
So $({\cal R}_M:
M\in\N)$ is a sequence of sets of density one. Hence,
$$
\eqalign{\pi(g) &=\lim\limits_{M\to\infty} {1\over M} \sum\limits_{m\in{\cal R}_M}
 {1 \over M} \P_{\underline w} \{ (\varphi^mx)_0 =g\}\cr
&=
\lim\limits_{M\to\infty} {1 \over M} \sum\limits_{m\in {\cal R}_M}
\P_{\underline w} \left\{ \sum\limits_{k\le m} {m \choose k} \mu^{m-k} 
\nu^k x_k=g \right\}.\cr}
$$
 From the Remark,  ${m \choose k}  \ne 0 \ mod \ p$ implies that
the mapping $\psi(g)= {m \choose k} \mu^{m-k}\nu^k\  g $ is one-to-one. Therefore from Lucas' theorem and
Lemma \equ(401) (a)  we get that for any $m \in {\cal R}_M$  
$$\left | \big\{k\le m: {m \choose k} \ mod \ p \ne 0 \big\} \right | 
\ge 2^{\alpha \log\log M}$$
and then 
$$
\left| \P_{\underline w}\left\{\sum\limits_{k\le m} {m \choose k} \mu^{m-k}\nu^k \ x_k=g\right\}-q^{-1}
\right| \le \vep_1(2^{\alpha\log\log M}).$$
Then $\pi(g) = q^{-1}$.
\medskip
Now we are ready to prove the result. Notice that for every 
$(g_j:j<s)\in G^s$ there exists a $(g'_j:j<s)\in G^s$ such that
$$
\{x\in G^\N:(\varphi^nx)_j=g_j\hbox{ for } j<s\}=\{x\in G^\N:
(\varphi^{n+j}x)_0=g'_j\hbox{ for }
j<s\} .
$$
Then it suffices to show that for any finite set $J\subseteq\N$ 
with $0\in J$ and $(g_j:j\in J)\in G^{J}$ it is verified,
$$
\lim\limits_{M\to\infty}{1\over M}\sum\limits_{m\le M}
\P_{\underline w}\{(\varphi^{m+j} x)_0=g_j, j\in J\} =q^{-|J|} .
$$
Introduce the following notation. We put 
$G_m=|\{n\le \delta_{1,m}:m_n<p-1\}|$ and we denote
$$
\{n\le \delta_{1,m}:m_n<p-1\}=\{\beta_{1,m}<\beta_{2,m}<...<\beta_{G_m,m}\}.
$$
Fix $\alpha\in(0,{1\over 2}),\varepsilon\in (0,\alpha),\varepsilon '\in 
(0,{1\over 2}(\alpha-\varepsilon))$. Denote $\ell =\max J$ and define
$$
{\cal R}'_M=\{ m\le M:\log(2(\ell+1))\le G_m\hbox{ and }
\beta_{[\log 2(\ell+1)],m}\le \varepsilon\log\log M\}.
$$
$$
{\cal R}''_M=\{ m\le M:\delta_{1,m}>\varepsilon\log\log M,\
|{\cal I}(m)\cap\{\varepsilon\log\log M \le n \le \delta_{1,m}\}|
\ge \varepsilon ' \log\log M\}.
$$
Both families of sets $({\cal R}'_M:M\in\N)$, $({\cal R}''_M:M\in \N)$ are 
of density 1.
\smallskip
Now for any family of sets $(\tilde {\cal R}_M:M\in \N)$ with 
$\tilde {\cal R}_M\subseteq \{m\le M\}$, we put $\tilde
{\cal R}_{M,J}=\{m\le M:
m+j\in \tilde{\cal R}_M$ for $j\in J\}$. If $(\tilde {\cal R}_M:M\in\N)$ 
is of density 1 then also
$(\tilde {\cal R}_{M,J}:M\in \N)$ is of density 1. Hence $({\cal R}_{M,J}:M\in \N)$, 
$({\cal R}'_{M,J}:M\in \N)$,
$({\cal R}''_{M,J}:M\in \N)$ are sequences of density 1.
\smallskip

Let $m\in {\cal R}'_{M,J}\cap {\cal R}''_{M,J}$. We denote 
${\cal I}_+(m+j)={\cal I}(m+j)\cap\{n>
\varepsilon \log\log M\}$, and ${\cal I}_-(m+j)={\cal I}(m+j)
\cap\{n\le \varepsilon\log\log M\}$. From the definition of ${\cal R}'_M$
we have that ${\cal I}_+(m+j)={\cal I}_+(m)$ for $j\in J$. Put
${\cal C}_+(m+j)=\{(m+j)_i:i\in {\cal I}_+(m+j)\}$ and
${\cal C}_-(m+j)=\{(m+j)_i:i\in {\cal I}_-(m+j)\}$ for $j\in J$. We have
${\cal C}_+(m+j)={\cal C}_+(m)$ for $j\in J$, and the sets
$({\cal C}_-(m+j):j\in J)$ are all different between them. Define for $j\in J$
$$
\tilde {\cal R}^j \! = \!
\{ 
k \le m+j: {\cal I}(k) \! \subseteq \! {\cal I}(m+j), k_i \! \le \! m_i
\hbox{ for } i \in \! {\cal I}_+(m),
k_i \! =  \! (m+j)_i \hbox{ for } i \in  \! {\cal I}_-(m+j) \}.
$$
The family $(\tilde {\cal R}^j:j\in J)$ is disjoint because the sets  
$({\cal C}_-(m+j):j\in J)$ are different. Moreover 
$|\tilde {\cal R}^j|\ge 2^{\varepsilon' \log\log M}$. 
\medskip
 From Lemma \equ(401) (b) and the Remark 
we get the result. In fact for every $m\in {\cal R}'_{M,J}\cap {\cal R}''_{M,J}$
and $j \in J$ we have that
$$
(\varphi^{m+j}x)_0=\sum\limits^{m+j}_{k=0} {m+j \choose k} \mu^{m+j-k}\nu^k x_k
$$
and the sequences 
$(\tilde {\cal R}^j:j\in J)$ satisfies conditions (H1),(H2),(H3). Indeed,
property (H1) follows from  
$\tilde {\cal R}^j \subset \{k\le m+j :{m+j \choose k} \ mod \ p   >0 \}
$,
they are disjoint, and if $k \in \tilde R^j$ then 
${m+j' \choose k} \ mod \ p=0$ for every $j'<j$ in  $J$ which shows property
(H3).
Then, from Lemma \equ(401) (b),  for any such $m$
$$
\left |\P_{\underline w}\{x:(\varphi^{m+j}x)_0=g_j, j\in J\}-q^{-|J|}\right |
\le \vep_J( 
2^{\varepsilon' \log\log M}).
$$
Then the theorem is shown. \square

\bigskip
\ni{\bf Acknowledgments.} Alejandro Maass and Servet Mart\'inez acknowledge 
financial support
 from C\'atedra Presidencial fellowship and Fondecyt grants 1980657 and
1970506.
Pablo A. Ferrari is partially supported by FAPESP
(Projeto Tem\'atico), CNPq (Bolsa de aux\'\i lio \`a pesquisa) and
FINEP (Projeto N\'ucleos de Excel\^encia).

\bigskip
\vfill\eject
\ni {\bf References.}
\bigskip
\item{[AN]} K.B. Athreya, P. Ney, {\it A new approach to the limit theory
of recurrent Markov chains}, Transactions of the AMS 248, 493--501 (1978).
\medskip

\item{[B]} H. Berbee,   {\it Chains with infinite connections: Uniqueness
and Markov representations}, 
Probab. Theory Related Fields 76 (1987), no. 2, 243--253.  
\medskip

\medskip

\item{[BFG]} X. Bressaud, R. Fern\'andez, A. Galves,  {\it Decay of
correlations for non H\"olderian dynamics. A coupling approach.} Preprint (1998).

\item{[FG]} P. Ferrari, A. Galves, 
{\it Acoplamento em processos estoc\'asticos}, 
21 Col\'oquio Bra\-si\-lei\-ro de Matem\'atica [21th Brazilian Mathematics Colloquium],
IMPA, Rio de Janeiro (1997). Available in
http://www.ime.usp.br/~pablo/abstracts/libro.html. 
\medskip

\item{[FMM]} P. Ferrari, A. Maass, S. Mart\'inez, 
{\it Ces\`aro mean distributionn of group automata 
starting from Markov measures}, 
Preprint (1998).
\medskip

\item{[L]} D. Lind, {\it Applications of ergodic theory and sofic systems 
to cellular automata}, Physica D 10, 36-44 (1984).
\medskip

\item{[Lu]} E. Lucas, {\it Sur les congruences des nombres eul\'eriens et des
coefficients diff\'erentiels des fonctions trigonom\'etriques,
suivant un module premier}, Bulletin de la Soc. Math\'ematique 
de France 6, 49--54 (1878).

\medskip

\item{[MM]} A. Maass, S. Mart\'inez, {\it On Ces\`aro limit distribution of
a class of permutative cellular automata}, Journal of Statistical Physics
90, 435--452 (1998).  

\item{[NN]} P. Ney, E. Nummelin, {\it Regeneration for chains with infinite
memory}, Probab. Theory Related Fields 96 (1993), no. 4, 503--520.

\bye


\magnification=1200
\def\'#1{{\if #1i{\accent"13\i}\else {\accent"13 #1}\fi}}

\hsize=32pc
\vsize=43pc
\vsize 19truecm
\baselineskip 14pt

\global\newcount\numsec\global\newcount\numfor

\def\senondefinito#1{\expandafter\ifx\csname#1\endcsname\relax}

\def\SIA #1,#2,#3 {\senondefinito{#1#2}%
\expandafter\xdef\csname #1#2\endcsname{#3}\else
\write16{???? ma #1,#2 e' gia' stato definito !!!!} \fi}

\def\etichetta(#1){(\veroparagrafo.\veraformula)%
\SIA e,#1,(\veroparagrafo.\veraformula) %
\global\advance\numfor by 1%
\write16{ EQ #1 ==> \equ(#1) }}

\def\letichetta(#1){\veroparagrafo.\veraformula
\SIA e,#1,{\veroparagrafo.\veraformula}
\global\advance\numfor by 1
 \write16{ Sta #1 ==> \equ(#1) }}

\def\tetichetta(#1){\veroparagrafo.\veraformula 
\SIA e,#1,{(\veroparagrafo.\veraformula)}
\global\advance\numfor by 1
 \write16{ tag #1 ==> \equ(#1)}}

\def\FU(#1)#2{\SIA fu,#1,#2 }

\def\etichettaa(#1){(A\veroparagrafo.\veraformula)%
\SIA e,#1,(A\veroparagrafo.\veraformula) %
\global\advance\numfor by 1%
\write16{ EQ #1 ==> \equ(#1) }}

\def\BOZZA{
\def\alato(##1){%
 {\rlap{\kern-\hsize\kern-1.4truecm{$\scriptstyle##1$}}}}%
\def\aolado(##1){%
 {
{
 \rlap{\kern-1.4truecm{$\scriptstyle##1$}}}}}
}

\def\alato(#1){}
\def\aolado(#1){}

\def\veroparagrafo{\number\numsec}\def\veraformula{\number\numfor}

\def\Eq(#1){\eqno{\etichetta(#1)\alato(#1)}}
\def\eq(#1){\etichetta(#1)\alato(#1)}
\def\leq(#1){\leqno{\aolado(#1)\etichetta(#1)}}
\def\teq(#1){\tag{\aolado(#1)\tetichetta(#1)\alato(#1)}}
\def\Eqa(#1){\eqno{\etichettaa(#1)\alato(#1)}}
\def\eqa(#1){\etichettaa(#1)\alato(#1)}
\def\eqv(#1){\senondefinito{fu#1}$\clubsuit$#1
\write16{#1 non e' (ancora) definito}%
\else\csname fu#1\endcsname\fi}
\def\equ(#1){\senondefinito{e#1}\eqv(#1)\else\csname e#1\endcsname\fi}

\def\Lemma(#1){\aolado(#1)Lemma \letichetta(#1)}%
\def\Theorem(#1){{\aolado(#1)Theorem \letichetta(#1)}}%
\def\Proposition(#1){\aolado(#1){Proposition \letichetta(#1)}}%
\def\Corollary(#1){{\aolado(#1)Corollary \letichetta(#1)}}%
\def\Remark(#1){\aolado(#1)Remark \letichetta(#1)}%

\def\include#1{
\openin13=#1.aux \ifeof13 \relax \else
\input #1.aux \closein13 \fi}

\openin14=\jobname.aux \ifeof14 \relax \else
\input \jobname.aux \closein14 \fi
\openout15=\jobname.aux

\def\include#1{
\openin13=#1.aux \ifeof13 \relax \else
\input #1.aux \closein13 \fi}
\openin14=\jobname.aux \ifeof14 \relax \else
\input \jobname.aux \closein14 \fi
\BOZZA

\def\one{{\bf 1}}
\def\ni{{\noindent }}
\def\sobre#1#2{\lower 1ex \hbox{ $#1 \atop #2 $ } }
\input amssym.def  
\input amssym.tex  
\def\E{{\Bbb E}}

\def\P{{\Bbb P}}

\def\vep{{\varepsilon}}
\def\Z{{\Bbb Z}}

\def\N{{\Bbb N}}
\def\bn{{\bf N}}
\def\bbn{{\bf \overline N}}
\def\R{{\Bbb R}}
\def\bm{{\bf M}}
\def\({\left(}
\def\){\right)}
\def\square{\ifmmode\sqr\else{$\sqr$}\fi}
\def\sqr{\vcenter{
         \hrule height.1mm
         \hbox{\vrule width.1mm height2.2mm\kern2.18mm\vrule width.1mm}
         \hrule height.1mm}}                  

\hfill{\bf Working Paper}
\bigskip
\centerline{\bf
 Ces\`aro 
mean distribution of group automata}
\centerline{\bf starting from measures with summable decay}
\vskip 3mm
\centerline{\bf Pablo A. Ferrari}
\centerline{\sl
Instituto de Matem\'atica e Estat\'istica,
Universidade de S\~ao Paulo}
\centerline{\sl
Caixa Postal 66281, 05315-970 S\~ao Paulo,
Brasil}
\centerline{\sl e-mail: \tt pablo@ime.usp.br}
\vskip 3mm
\centerline{\bf Alejandro Maass, Servet Mart\'inez}
\centerline{\sl
Departamento de Ingenier\'{\i}a Matem\'atica}
\centerline{\sl Universidad de Chile,
Facultad de Ciencias F\'{\i}sicas y
Matem\'aticas,}
\centerline{\sl
Casilla 170-3, Correo 3, Santiago, Chile}
\centerline{\sl e-mail: \tt 
amaass@dim.uchile.cl, smartine@dim.uchile.cl}
\vskip 3mm
\centerline{\bf Peter Ney}
\centerline{\sl
Department of Mathematics}
\centerline{\sl University of Wisconsin,}
\centerline{\sl
Madison, WI 53706, U.S.A.}
\centerline{\sl e-mail: \tt 
ney@math.wisc.edu}
\vskip 3mm
\vskip 3mm
\centerline{\bf Abstract}
\bigskip
\bigskip

{\it Consider a finite Abelian group $(G,+)$, with $|G|=p^r$, $p$ a
prime number, and $\varphi: G^\N \to
G^\N$ the cellular automaton given by $(\varphi x)_n=\mu x_n+\nu
x_{n+1}$ for any $n\in \N$, where $\mu$ and $\nu$ are integers
relatively primes to $p$.  We prove that if $\P$ is a translation
invariant probability measure on $G^\Z$ determining a chain with
complete connections and summable decay of correlations, then for any
${\underline w}= (w_i:i<0)$ the Ces\`aro mean distribution
$\displaystyle {\cal M}_{\P_{\underline w}} =\lim_{M\to\infty} {1\over
M} \sum^{M-1}_{m=0}\P_{\underline w}\circ\varphi^{-m}$, where
$\P_{\underline w}$ is the measure induced by $\P$ on $G^\N$
conditioning to $\underline w$, exists and satisfies ${\cal
M}_{\P_{\underline w}}=\lambda^\N$, the uniform product measure on
$G^\N$. The proof uses a regeneration representation of $\P$ in the
vein of Ney and Nummelin.}

\bigskip
\ni{\bf 1.-- Introduction and main results.}
\numsec=1\numfor=1\bigskip
\bigskip

Let $(G,+)$ be a finite Abelian group with 
$q=p^r$ elements, being $p$ a prime number.  We put 
$\lambda=(q^{-1},..., q^{-1})$ the uniform measure on the
group.  In this paper we study the measure evolution under the
dynamics of 
the cellular automaton $\varphi:G^\N\to G^\N$, given by
$(\varphi x)_n=\mu x_n+ \nu x_{n+1}$ for $n\in\N$, where $\mu$ and
$\nu$ are integers relatively primes to
 $p$ ($\ell g $ means $g+...+g$ $\ell$--times).  
The uniform product measure
$\P=\lambda^\N$ is $\varphi$--invariant,  
$\P\circ\varphi^{-n}=\P$, but any other product measure $\P=\pi^\N$,
with $\pi\neq \lambda$, is not $\varphi$-invariant. Moreover, even in
the simplest case $G=\{0,1\}$ and $+$ the mod 2 sum, the limit of the
marginal distribution, $\lim\limits_{m\to\infty}\P \{ (\varphi^mx)_0=g
\}$ with $g \in G$, does not exist. 
The reason is that for $m=2^k$,
$(\varphi^mx)_0 = x_0+x_m$ (the other terms sum an even number of
times and do not contribute to the sum) has probability
$p^2+(1-p)^2$ to be $0$,
while for $m=2^{k}-1$ this probability converges to $1 \over 2$
because $(\varphi^mx)_0=\displaystyle\sum_{\ell=0}^m x_\ell$.

Alternatively we can study the Ces\`aro mean distribution
$$
{\cal M}_\P\,\doteq\,\lim_{M\to\infty} {1\over M} \sum^{M-1}_{m=0}
\P\circ\varphi^{-m}
$$ 
for a class of initial distributions $\P$ on $G^\N$. In the above
display and in the sequel $\doteq$ means ``it is defined by''.

 Let $-\N^* =\{-i:i\in\N\setminus\{0\}\}$ and $\N^*=\N\setminus\{0\}$.
Let $\P$ be a translation invariant probability measure on $G^\Z$. 
For $\underline w \in G^{-\N^*}$ let $\P_{\underline w}$
be the measure on $G^\N$ induced by the conditional probabilities as 
follows. For
any $m\ge 0$ and $g_0,\dots,g_m\in G$, define
$$
\P_{\underline w}(x_0=g_0,\dots,x_m=g_m)\,
\doteq\, \P(x_0=g_0,\dots,x_m=g_m\,\vert\,
x_{i}=w_{i}, i<0).
$$
We
say that $\P$ has {\sl complete connections} if it satisfies
$$
\hbox{ $\forall g_0 \in G, \ \forall \underline w \in G^{-\N^*}$, 
\quad $\P_{\underline w}\{x_0=g_0\}>0. $ }\Eq(50)
$$
For any $m\ge 0$ define
$$ 
\gamma_m\, \doteq\, \sup\left\{\left |{\P_{\underline w}\{ x_0=g
\} \over \P_{\underline v}\{ x_0=g \}} -1 \right | : g \in G, \underline v,
\underline w \in G^{-\N^*}, 
v_{i}=w_{i},i\in[-m,-1]\right\}.
$$ 
We say that $\P$ has {\sl summable decay} if
$$
\sum_{m=0}^{\infty} \gamma_m<\infty. \Eq(60)
$$
This is a uniform continuity condition on $\P_{\underline w}(g)$ as a function
of $\underline w$.

\medskip

The Ces\`aro limits has been already studied for the mod 2 sum
automaton and other classes of permutative cellular automata in [L]
and [MM].  In these papers it is computed mainly for Bernoulli
measures, and in [MM] only the one site Ces\`aro limit is computed for
a Markov measure.  In the mod 2 case the limit is uniformly
distributed, but for some permutative cellular automata the Ces\`aro
mean exists but it is not necessarily uniform. In [FMM] the
Athreya-Ney regeneration times representation of $r$-step Markov chain
was used to show the convergence of the Ces\`aro mean of the group
automata starting with these Markov chains to the uniform Bernoulli
measure.

In this paper we generalize these results for the group automaton
$\varphi$ and initial measures with complete connections and summable
decay.  

\bigskip
\proclaim \Theorem (1).
Let $(G,+)$ be a finite Abelian group with $|G|=p^r$, $p$ a prime
number and $r \in \N\setminus \{0\} $.  Let $\P$ be a translation 
invariant measure on
$G^\Z$ with complete connections and summable decay.  
Let 
$\varphi:G^\N\to G^\N$ be the cellular automaton, given by
$(\varphi x)_n=\mu x_n+ \nu x_{n+1}$ for $n\in\N$, where $\mu$ and
$\nu$ are integers relatively primes to $p$.  
Then for all
$\underline w \in G^{-\N^*}$ the Ces\`aro mean distribution ${\cal
M}_{\P_{\underline w}}$ exists and verifies ${\cal M}_{\P_{\underline
w}}= \lambda^\N$, the product of uniform measures on $G$. 

\bigskip 

There are two main elements in the proof: regeneration times
and distribution of Pascal triangle coefficients mod $p$.

\bigskip
\noindent{\bf 2.-- Regeneration times for the initial measure.} 
\numsec=2\numfor=1
\bigskip

We show that under the conditions of Theorem \equ(1), for all $\underline w\in
-\N^*$ we can jointly construct a random sequence $\underline x=(x_i:i \in \N)
\in G^\N$ with distribution $\P_{\underline w}$ and a random subsequence
$(T_i:i\in\N^*)\subseteq \N$ such that $(x_{T_i}:i\in \N^*)$ are iid uniformly
distributed in $G$ and independent of $(x_n:n\in \N^*\setminus
\{T_1,T_2,\dots\})$; furthermore $(T_i:i\in\N^*)$ is a stationary renewal
process with finite mean inter-renewal time independent of $\underline w$. A
consequence of the construction is that the random vectors (of random lenghts)
$((x_{T_i},\dots, x_{T_{i+1}-1}): i\ge 1)$ are iid. 

Our regeneration approach shares results with Berbee (1987) and Ney and
Nummelin (1993). The construction is simple: the probability space is
generated by product of iid uniform (in $[0,1]$) random variables. It works as
the well known construction and simulation of Markov chains as a function
of a sequence of uniform random variables (see, for instance Ferrari and
Galves (1997)). Bressaud, Fern\'andez and Galves (1998) construct a coupling
using these ideas to show decay of correlations for measures with infinite
memory.

For $\underline w \in G^{-\N^*}$ and $g \in G$ denote
$$
P(g\vert \underline w)\doteq \P\{x_0=g\vert x_i=w_i,
i\le -1\}.
$$

Let  
$$
a_{-1}(g\vert \underline w)\,\doteq \,
\inf \{ P(z\vert \underline v)\ : \  {\underline v} \in G^{-\N^*},\, z \in 
G \}. \Eq(g1)
$$ Actually $a_{-1}$ depends neither on $g$ nor on $\underline w$; we keep the
dependence in the notation for future (notational) convenience.  Since the
space $G^{-\N^*}$ is compact, the infimun in \equ(g1) must be attained by a
$g^0\in G$ and a $\underline w^0\in G^{-\N^*}$. Hence,
$$
a_{-1}(g\vert \underline w)\, =\,P(g^0\vert \underline w^0) \,>\, 0,
$$
because $\P$ has complete connections.
For each $k \in\N$, $g \in G$ and ${\underline w} \in G^{-\N^*}$
define 
$$
a_k(g \vert {\underline w}) \doteq 
\inf\{P(g\vert w_{-1},\dots, w_{-k}, {\underline z}) \ : \ 
{\underline z} \in G^{-\N^*} \},
$$
where $(w_{-1},\dots, w_{-k}, {\underline z}) = (w_{-1},\dots, w_{-k},
z_{-1}, z_{-2},\dots)$. Notice that $a_0(g\vert \underline w)$ does not depend
on $\underline w$. 
Let 
$$
b_{-1}(g\vert \underline w) 
\,\doteq\, a_{-1}(g\vert \underline w) 
,
$$ 
for $g\in G$. For $k\ge 0$,
$$ 
b_k(g \vert {\underline w}) \,\doteq\, a_k(g \vert {\underline
w})-a_{k-1}(g \vert {\underline w}).
$$

We construct disjoint intervals $B_k(g \vert {\underline w})$ for
$g\in G$, $k\ge -1$, contained in $[0,1]$, of Lebesgue measure $b_k(g
\vert {\underline w})$ respectively, disposed in increasing order with
respect to $g$ and $k$: $ B_{-1}(0|\underline w), \dots,B_{-1}(q-1|\underline w),
B_0(0|\underline w),\dots,B_0(q-1|\underline w),B_1(0|\underline
w),\dots, B_1(q-1|\underline w), \dots$, with no intersections (we
have enumerated $G$ by $\{0,...,q-1\}$). The construction guarantees
$$
\Bigl|\bigcup_{k\ge -1} B_k(g \vert {\underline
w})\Bigr|\,=\,P(g|\underline w)
$$
and
$$ 
\Bigl|\bigcup_{g\in G}\,\bigcup_{k\ge -1} B_k(g
\vert {\underline w})\Bigr|\,=\, 1.
$$ 
(All the unions above are disjoint.)

Let $\underline U= (U_n: n \in \Z)$ be a double infinite sequence of
iid random variables uniformly distributed in $[0,1]$. Let
$(\Omega,{\cal F},\P)$ be the probability space induced by these
random variables. For each $\underline w\in G^{-\N^*}$ we construct
the random sequence $\underline x$ with distribution $\P_{\underline
w}$ in $\Omega$, as a function of $\underline U$, recursively: for
$n\in \N$
$$ 
x_{n} \,\doteq \,\sum_{g\in G} g \Bigl[ \sum_{\ell\ge -1} \one\{ U_n \in
B_\ell(g \vert x_{n-1},\dots,x_0,{\underline w}) \}\Bigr].
$$ 
where we adopted the convention $B_\ell(g \vert
x_{n-1},\dots,x_0,{\underline w}) = B_\ell(g \vert
{\underline w})$. 

For $\ell\ge -1$ let
$$
B_\ell(\underline w) \doteq 
\bigcup_{g\in G}  B_\ell(g|\underline w).
$$
Notice that neither $B_{-1}(g|\underline w)$ nor $B_{-1}(\underline
w)$ depend on $\underline w$. Furthermore
$$
{|B_{-1}(g|\underline w)|\over|B_{-1}(\underline w)|} \, =\, |G|^{-1}.\Eq(gm1)
$$
For $k \in \N$ let
$$
a_k \doteq \min_{\underline w} \ \left \{ \displaystyle  \sum_{g\in G} 
 a_k(g \vert \underline w)\right\}.
$$
This is a non-decreasing sequence and satisfies
$$
[0,a_k]\;\subset \;\bigcup_{\ell=-1}^k 
B_\ell(\underline w), \Eq(oak)
$$ 
independently of $\underline w \in G^{-\N^*}$.
\bigskip
\proclaim\Lemma(71). In the event $\{U_n\le a_k\}$ for $n\in \N$ 
we only need to look at
$x_{n-1},\dots,x_{n-k}$ to decide the value of $x_n$. More precisely,
for $\underline v\in G^\Z$ such that $v_i=w_i$ for $i\le -1$,
$$
\eqalign{
\P_{\underline w}\{x_n=g \,\vert\, &U_n\le a_k,\,x_{n-1}=v_{n-1}, \dots, 
x_0=v_0\}\,\cr
&= \,\P_{\underline w}\{x_n=g
\,\vert\, U_n\le a_k,\,x_{n-1}=v_{n-1}, \dots, x_{n-k}= v_{n-k}\}.}
$$

\noindent {\bf Proof.} Follows from \equ(oak). \square
\medskip

Define times
$$
\eqalign{
T_1&\doteq \min \{n\ge 0: \,U_{n+j} \le a_{j-1}, \
j\ge 0 \},\cr
T_i &\doteq \min \{n>T_{i-1}: \,U_{n+j} \le a_{j-1}, \
j\ge 0 \}, }
$$ 
for $i > 1$.  
\bigskip
\proclaim\Lemma(72). The variables $(x_{T_i}:\,i\ge 0)$ are iid
uniformly distributed in $G$.

\noindent {\bf Proof.} Let us show that the marginal distribution of
$x_{T_i}$ is uniform in $G$:
$$
\eqalign{
\P\{x_{T_i}=g \}
&= \sum_{n\in \N}\P\Bigl \{U_n \in \displaystyle \bigcup_{\ell \ge -1}
B_\ell(g\vert \underline w) \,,\,T_i = n\Bigr \}\cr
&= \sum_{n\in \N}\P \{U_n \in B_{-1}(g\vert\underline w)\,\vert\, U_n\in
B_{-1}(\underline w)\}\, \P\{T_i = n\} \cr
&= |G|^{-1}.}  
$$ 
The second identity follows because $\{T_i=n\}$ is the intersection
of $\{U_n \in  B_{-1}(\underline w)\}$ with events depending on variables
$(U_{n+\ell}, \, \ell\neq 0)$ which are independent of $U_n$. The
third identity follows from \equ(gm1). The
same computation shows that for any $K\subset \N$, $(i(k):k\in
K)\subseteq \N$, and $(g_k:k\in K)\subseteq G^{|K|}$
$$
\P\{x_{T_{i(k)}} = g_k,\, k\in K\} = |G|^{-|K|},
$$ 
so that $(x_{T_{i(k)}}:k\in K)$ are iid in $G$.  The reason why
the above computation works is that in the event $\{T_i=n\}$,
$U_{n+1}\le a_0$, hence $x_{n+1}$ does not depend on the past. Since for
all $j\ge 1$, $U_{n+j}\le a_{n+j-1}$, $x_{n+j+1}$ only depends on
$x_{n+1},...,x_{n+j}$.
\square

\medskip

Let $\bn$ be the counting measure on $\N$ induced by $(T_i:i\ge 1)$:
for $A\subset \N$ and $n\in\N$, 
$$ \bn(A) \,\doteq\, \sum_{i\ge 1} \one\{T_i\in A\},\ \ \
\bn(n)\,\doteq\,\bn(\{n\}).
$$
Notice that the definitions of $(T_i:i\ge 1)$ and $\bn$ depend only on $(U_n:n\in \Z)$
and do not depend on~$\underline w$.

\proclaim\Lemma(70). The distribution of the counting measure $\bn$
corresponds to a stationary renewal process.

\noindent{\bf Proof. } We will construct a stationary renewal process
$\bm$ in $\Z$ whose projection on $\N$ is $\bn$.
For $k\in\Z$, $k'\in\Z\cup\{\infty\}$, define
$$
H[k,k'] \,\doteq\,\cases{  \{U_{k+\ell}\le a_{\ell-1},
\ell=0,\dots,k'-k\},&if $k\le k'$\cr
&\cr
\hbox{``full event''},&if $k>k'$ }
$$
With this notation,
$$
\bn(n) = \one \{H[n,\infty]\},\ \ n\in\N. \Eq(500)
$$
We construct a double infinity counting process $\bm$ using the
variables $(U_n:n\in \Z)$ by 
$$
\bm(n) \,\doteq\, \one \{H[n,\infty]\},\ \ n\in\Z.
$$ 
By construction, the distribution of $\bm$ is translation
invariant, hence $\bm$ is stationary. Furthermore, by \equ(500) it
coincides with $\bn$ in $\N$: $\bm(K)=\bn(K)$ for $K\subset \N$. 
Define $T_i$ for $i\le 0$ as the ordered time-events of $\bm$ in the
negative axis. 

The (marginal) probability of a counting event at time $n\in\Z$, is
given by
$$ 
\P\{\bm(n)=1\}
= \P\{U_{n+j} 
\le a_{j-1}, \ j\ge 0\}
\,=\, a_{-1}\, a_0\, a_1\,\cdots 
\,\doteq\, \beta,
$$
and it is independent of $n$. 
We first show that under the hypothesis of summability of $\gamma_k$,
$\beta $ is strictly positive. 
For any $g \in G$, $w_{-1},\dots, w_{-k} \in G$ 
and ${\underline z}, {\underline v}
\in  G^{-\N^*}$
$$
\left | {\P\{g\vert  w_{-1}\dots w_{-k}, {\underline z}\} \over 
\P\{g\vert  w_{-1}\dots w_{-k}, {\underline v}\}}
-1 \right | \le \gamma_k,$$
therefore
$$
\inf \ \{ \P\{g \vert  w_{-1}\dots w_{-k}{\underline z}\} \ : \
{\underline z} \in G^{-\N^*} \}\ge 
(1-\gamma_k) \P\{g\vert  w_{-1}\dots w_{-k} {\underline v}\}.
$$
Summing over $g\in G$ and taking minimum on the set
$\{w_{-1},\dots,w_{-k}\}$ we conclude that
$$a_k \ge 1-\gamma_k.$$ 
Since $\displaystyle 
\sum_{k\ge 0} \gamma_k < \infty$ we deduce that 
$\displaystyle 
\sum_{k\ge 0} (1-a_k) < \infty$ and henceforth 
$\beta >0$.

\medskip

We show now that $\bm$ is a renewal process on $\Z$.
The event $\{\bm(n) = 1\}$ depends only on $(U_k:k\ge n)$, that is,
$(T_i:i\in\Z)$ are stopping times for the process $(U_{-k}:k\in \Z)$.
Since for $k<k'<k''\le \infty$, 
$$
H[k,k'']\cap H[k',k'']\,\, =\,\, H[k,k'-1]\cap H[k',k''],
$$ we have that for any finite set $A=\{k_1,\dots,k_n\}$ with
$k_1<\dots<k_n<k'$ and for any sequence $(m_\ell:\ell>k')$ with
$m_\ell\in\{0,1\}$,
\smallskip
$$
\eqalign{
\P\big\{&\bm(A)=n\,\big\vert\, \bm(k')=1, \bm(\ell)= m_\ell,
\,\ell>k'\big\}\cr
&=\,\P\big\{\displaystyle \bigcap_{i=1}^nH[k_i,k'-1]\,\big\vert\, \bm(k')=1\big\}
\phantom{\sum_o^p}\cr 
&= \,\prod_{i=1}^n\P\{H[k_i,k_{i+1}-1]\},\cr
}\Eq(epb)
$$ 
where $k_{n+1}\doteq k'$. The computation above could be done
because $\P\{\bm(k')=1\}=\beta>0$.  Display \equ(epb) means that given a
counting event at time $k'$, the distribution of the counting events
for times less than $k'$ does not depend on the events after
$k'$. This characterizes $\bm$ as a renewal process.  Since the
density $\beta$ is positive, $T_1$, the residual time is a honest
random variable, and for $i\neq1$,
$\E(T_{i+1}-T_i)=\beta^{-1}<\infty$.\ \ \square

\medskip\medskip
\ni{\bf 3.-- A renewal Lemma.}
\medskip\medskip
\numsec=3 \numfor=1

In this section we show that a stationary discrete-time renewal
process on $\N$ has high probability to visit sets with many points.

\proclaim \Lemma (4).  Let $\bn$ be a stationary renewal process with
finite inter-renewal mean. Then for all $A\subset \N$,
$$
\P\{\bn(A)=0\} \le \vep(|A|)
$$ 
with $\vep(n)\to 0$ as $n\to\infty$. Also, $\vep :\N \to \R$ can be
chosen to be decreasing.


\noindent{\bf Proof.} We are going to prove that for all $\vep>0$
there exists $n_0$ such that for any finite set $A\subset \N$ with
$|A|>n_0$,
$$
\P\{\bn(A)=0\} \le \vep. \Eq(221) 
$$ 
We start with some known facts of renewal theory.  Let $T_i$ be the
renewal times and $\beta = 1/\E(T_{i+1}-T_i)$ for some $i\ge 1$ (and
hence for all $i\ge 1$). Since the inter-renewal distribution has a
first moment finite, the key renewal theorem holds: we have
$$
\lim_{n\to\infty}\P\{\bn(n) =1\,\vert\, \bn(0)=1\}\,=\,  \beta. \Eq(220)
$$
Let $S_n\doteq T_{\bbn(n)+1}-n$ be the residual time (over jump) at $n$,
where we have denoted by $\bbn(n)=\sum_{k\le n}\bn(k)$,
and let for $k\ge 0$
$$
\beta_k=\P(\bn(k)=1\,|\, T_1=0),\ \ \ F(k) = \P(T_2-T_1>k),\ \ \ 
F_n(k) = \P(S_n >k). \Eq(221)
$$
Now we have
$$
F_n(k)\,=\, \sum_{j=0}^n F(j+k) \beta_{n-j}\,\le\, \overline F(k), \Eq(222)
$$
where
$$
\overline F(k)\,\doteq\, \sum_{j=k}^\infty F(j)\,\to \,0
$$
as $k\to\infty$ because we are assuming that the inter-renewal time
has a finite mean.

For any subset $B\subset A$ we have
$$
\P\{\bn(A)=0\}\,\le\,\P\{\bn(B)=0\}.
$$
For any $A$ with $|A|=n$ and any $1<\ell<n$, there exists a set
$$
\{b^n_1,\dots,b^n_\ell\} \doteq B^n_\ell\subset A
$$
with
$$
\Bigl[{n\over \ell}\Bigr]\,\le\, b^n_{j+1}-b^n_j, \ \ \ j=1,\dots,\ell-1,\Eq(224)
$$
where $[x]$ is the largest integer in $x$. The choice of 
$\{b^n_1,...,b^n_\ell\}$
depends on $A$ but $\ell$, and \equ(224) hold uniformly for all $A$
with $|A|=n$.

Let $\vep>0$ and take any $0<\delta<\beta$. Take $n_0$ such that
$\beta_n>\delta$ for $n>n_0$. Let  $n>\ell n_0$ and define 
$$ 
\Gamma^n_j \,\doteq\, \{S_{b^n_j}\,\le [n/\ell] - n_0\}
$$
the event ``the over jump of $b^n_j$ does not superate $[n/\ell]-n_0$''.
Let
$$
\Theta^n_j\,\doteq\,
\bigl\{\bn(b^n_j-b^n_{j-1}-S_{b^n_{j-1}})\,=\,0\bigr\}
$$
the event ``starting at the over jump of $b^n_{j-1}$, $b_j^n$ is not
hit''.

Now
$$ 
\eqalign{ \P\{\displaystyle\bigcap_{j=1}^\ell \Gamma^n_j \}\,
&=\, \P\{\hbox{``all overjumps of $b^n_j$'s are $\le [n/\ell]-n_0$''}\}\cr  
&=\,\P\{\Gamma^n_\ell\,|\,\displaystyle\bigcap_{j=1}^{\ell-1}
\Gamma^n_j\}\,\P\{\displaystyle\bigcap_{j=1}^{\ell-1} \Gamma^n_j\} \cr 
&=\,\P\{\Gamma^n_\ell\,|\,\Gamma^n_{\ell-1}\}\,
\P\{\Gamma^n_{\ell-1}\,|\,\Gamma^n_{\ell-2}\}\dots
\P\{\Gamma^n_2\,|\,\Gamma^n_{1}\}\,\P\{\Gamma^n_{1}\}\cr
&\ge \big(1- \overline F([n/\ell]-n_0)\big)^\ell}\Eq(226)
$$
by \equ(222). Then
$$
\eqalign{
\P\{\bn(A) =0\} \,&\le\, \P\{\bn(B^n_\ell)=0\}\cr
&\le \,\P\{\{\bn(B^n_\ell)=0\}\,\cap\,\displaystyle\bigcap_{j=1}^{\ell} 
\Gamma^n_j\}
\,+\,
1-\P\{\displaystyle\bigcap_{j=1}^{\ell} \Gamma^n_j\}.} 
$$
Now
$$
\eqalign{
\P\{\{\bn(B^n_\ell)=0\}\,\cap\,\displaystyle\bigcap_{j=1}^\ell \Gamma^n_j\} 
&= \P\{\displaystyle\bigcap_{j=1}^\ell (\Gamma^n_j\cap\Theta^n_j)\}\cr
& \le \P\{\displaystyle\bigcap_{j=1}^\ell \Theta^n_j \, \big | \,
\bigcap_{j=1}^\ell \Gamma^n_j \}  
\P\{\displaystyle\bigcap_{j=1}^\ell \Gamma^n_j\} \cr
& \le \P\{\displaystyle\bigcap_{j=1}^\ell \Theta^n_j\}
 \le\, (1-\delta)^\ell. \cr }\Eq(225)
$$
since $\beta_n > \delta$ for $n > n_0$.
By \equ(226) and \equ(225),
$$ 
\P\{\bn(B^n_\ell)=0\}\,\le\,(1-\delta)^\ell\,+\, 1 -
\bigl(1-\overline F([n/\ell]-n_0)\bigr)^\ell.
$$
Now choose $\ell$ so that $(1-\delta)^\ell<\vep/2$, then $n$ so that
$1-(1-\overline F([n/\ell]-n_0))^\ell<\vep/2$, to conclude
$$
\P\{\bn(B^n_\ell)=0\}\,\le\,\vep
$$
for sufficiently large $n+\ell$. \square

\bigskip
\bigskip
\ni{\bf 4.-- Convergence of Ces\`aro limit.}
\bigskip
\bigskip
\numsec=4 \numfor=1

For proving this theorem we shall need some results concerning 
walks of variables
determining a chain with complete connections. 
In this purpose let 
us introduce some notation.
First $R=(r_k:k\in\N)$ denotes an increasing sequence in $\N$. 
We put  
$R_n=(r_k:k\le n)$. For any subsequence 
$\overline R=(\overline r_k: k \in \N)$ of $R$ we define the index function by
$f_{\overline R}(k)=\ell $ if $\overline r_k =r_\ell$. We also set 
$n(\overline R)=|\overline R\cap R_n|$.
Let $a^R=(a^R_r:r\in R)$ be a sequence of non-negative integers. They
define maps  $\psi^R_r:G\to G$ such that $\psi_r^R(g)=a_r^R \ g=g+...+g  $ \
$a_r^R$ times, for any $r \in R$.
We associate to it the following 
sequence of random variables taking values in $G$,
$$
S^R_n =\sum\limits_{r\in R_n} a_r^R \ x_r,\ n\in\N. 
$$
We will distinguish  the following subsequence 
$$R^* \doteq R^*(a ^R) = \{ r\in R \ : \ a_r^R\neq 0 \ mod \ p  \}.$$

\medskip
\ni {\bf Remark.} Since $(G,+)$ is a finite Abelian group with 
$|G|=p^r$, $p$ a prime number, then the function $\psi(g)=a\ g $, 
where  $a\in \N$, is one-to-one whenever $a\neq 0 \ mod \ p$.
\medskip

Let $J\subseteq \N$ be a finite set. Consider 
a finite family of sequences $R^J=(R^j:j\in J)$. Associated to each
sequence there is a sequence of non-negative 
integers $a^{R^j}=(a_r^{R^j}:r \in R^j)$ and the corresponding
set of mappings
$\psi^{R^j}= (\psi^{R^j}_r:r\in R^j)$.
As before we consider the sequences
$R^{j*}\doteq R^*(a^{R^j})$ for 
$j\in J$. Let $\tilde R^J=(\tilde R^j:j\in J)$ 
be a family of subsequences verifying the following conditions:
\medskip
\item{(H1)} $\tilde R^j\subseteq R^{j*}$ for any $j \in J$, 
\item{(H2)} $\tilde R^j \cap \tilde R^i=\emptyset$ 
if $i \neq j$ in $J$, 
\item{(H3)} if $r \in \tilde R^j\cap R^k $ for $k<j$ in $J$, then
$a_r^{R^k}=0 \ mod \ p$. 
\medskip
\ni We set 
$$\tilde n(\tilde R^J)=\min\{ n(\tilde R^j):j\in J\}$$
and
$$\tilde n(R^J)=\max\{\tilde n (\tilde R^J):\tilde R^J 
\hbox{ verifying } (H1), (H2), (H3) \}.$$
\medskip
The proof of  Theorem \equ(1) is based upon the   following result.

\bigskip
\proclaim \Lemma (401). 
Let $\P$ be a translation invariant measure on $G^\Z$ 
with complete connections such that
$\displaystyle\sum_{m\ge 0} \gamma_m < \infty$, and let $\underline w \in G^{-\N^*}$.
Then
\item{(a)} $\exists \ \vep_1:\N \to \R$, a decreasing function with 
$\vep_1(n)\to 0$ if $n \to \infty$,  
such that  for any increasing sequence $R$ in $\N$ and any
sequence of non-negative integers $a^R$ it is verified 
$$
\left |\P_{\underline w}\{ S^R_n=g\}-q^{-1} \right | \le \vep_1(n(R^*)), \ \hbox{ for any } 
 n\in \N, g\in G.
$$
\item{(b)} Let $J \subset \N$ be  finite.  Then there is a decreasing function
$\vep_J:\N \to \R$ with  $\vep_J(n)\to 0$ if $n\to \infty$, 
such that  
for any set of sequences $R^J=(R^j:j\in J)$ and any  family of
non-negative integers   
$(a^{R^j}:j\in J)$, it is verified 
$$\left |\P_{\underline w}\{S^{R^j}_n=g_j, \hbox{ for }  j\in J\} -q^{-|J|}\right |\le 
\vep_J(\tilde n(R^J)) 
\hbox{ for any }n\in \N, (g_j:j\in J)\in G^{J} .
$$
\square

\bigskip
Before begin the proof of Lemma \equ(401) we include a useful arithmetic property. We include
a proof for completeness.
For $(G,+)$ a finite  Abelian group with $|G|=p^r$, where
$p$ is a prime number, consider the following system 
of equations (S):
\bigskip
$$\matrix{
(1)&a_{11}g_1 &+&a_{12}g_2&+&...&+&a_{1\ell}g_{\ell}=0 \cr
(2)&a_{21}g_1 &+&a_{22}g_2&+&...&+&a_{2\ell}g_{\ell}=0 \cr
   &          & &         & &\vdots   & &                    \cr
(\ell)&a_{\ell1}g_1 &+&a_{\ell 2}g_2&+&...&+&a_{\ell \ell}g_{\ell}=0 \cr
}
$$
such that 
$$
(H') \qquad a_{ij} \in \N, \ a_{ii}\neq 0 \ mod \ p, \ a_{ij}=0 \
mod \ p  \hbox{ if } i < j.
$$ 
Denote $a_{ii}=k_ip+s_i$ with $s_i \in \{1,...,p-1\}$
and $a_{ij}=c_{ij}p$ for $i<j$.
\bigskip

\proclaim \Lemma (402).  
The system (S) has unique solution $g_1=g_2=...=g_{\ell}=0$.

\smallskip
\ni {\bf Proof.}
First of
all we will prove that if $g_1,...,g_{\ell}$ are solutions of
(S) and for some $1 < s \le r$,  $p^sg_i=0$, $i \in \{1,...,\ell\}$,
then $p^{s-1}g_i=0$ for $i \in \{1,...,\ell\}$. 
We prove this property by induction on $\{1,...,\ell\}$.
First consider equation (1),
$$(k_1p+s_1) g_1+\sum_{j=2}^{\ell}c_{1j}pg_j=0.$$
If we add the equation $p^{s-1}$ times we obtain,
$$k_1p^sg_1+s_1p^{s-1}g_1+\sum_{j=2}^{\ell}c_{1j}p^sg_j=0,$$
then $s_1p^{s-1}g_1=0$. Since the product by $s_1$ defines 
a 1-to-1 map we conclude that $p^{s-1}g_1=0$.
Let us  continue with the induction assuming that $p^{s-1}g_1=0$,
$p^{s-1}g_2=0,...,p^{s-1}g_t=0$, for $1\le t < \ell$, and we prove that
$p^{s-1}g_{t+1}=0$.
\medskip
Adding $p^{s-1}$ times equation $t+1$ we get
$$\sum_{j=1}^ta_{t+1,j}p^{s-1}g_j+(k_{t+1}p+s_{t+1})p^{s-1}g_{t+1}+
\sum_{j=t+2}^{\ell}c_{t+1,j}p^{s}g_j=0.$$
Therefore, using the induction hypothesis we obtain
$s_{t+1}(p^{s-1}g_{t+1})=0$ and henceforth $p^{s-1}g_{t+1}=0$.
\medskip
To conclude we use last property recursively
beginning from the fact that $p^rg_i=0$ for
any $i \in \{1,...,\ell\}$. \square

Hence the transformation $A:G^\ell \to G^\ell,$ $A \vec g=\vec h$, with $\vec g,
\vec h \in G^\ell$ and matrix $A$ verifying condition (H') is a one-to-one
and onto transformation. In what follows we identify $\P_{\underline w}$
with $\P$.
\bigskip
\ni {\bf Proof of Lemma \equ(401).}
\bigskip

\noindent a) For any increasing sequence $R= (r_k:k\in\N)$ we put
$$
\tau^{R} = \inf\{k\in \N: \bn(r_k)=1\}, 
\hbox{ where } \infty=\inf\phi, 
$$
the first time that some element of the sequence $R$ belongs to the
renewal process $\bn$ introduced in Section 2. Consider $R^*$ the
subsequence corresponding to one-to-one mappings associated to $R$. We
denote $n^*=n(R^*)$, $\tau^*=\tau^{R^*}$ and $f=f_{R^*}$ the
corresponding index function.  First we prove
$$
\P\{S_n^R=g\vert \tau^{*}\le n^*\}=q^{-1} .
$$
To see that write
$$
\eqalign{
&\P \left \{ S_n^R=g, \tau^{*}\le n^* \right \}= \sum_{k=1}^{n^*} 
\P\left \{S_n^R=g, \tau^{*}=k \right \}\cr
&=\sum_{k=1}^{n^*} \P \left\{\sum_{i=1}^{f(k)-1}\psi_{r_i}(x_{r_i})+
\psi_{r_{f(k)}}
(U_{r_{f(k)}})+
\sum_{i=f(k)+1}^{n}\psi_{r_i}(x_{r_i})=g, \tau^{*}=k\right\} \cr
&=\sum_{k=1}^{n^*}\sum_{g_1, g_2 \in G} 
\P \Bigl \{  \sum_{i=1}^{f(k)-1}\psi_{r_i}(x_{r_i})=g_1, 
U_{r_{f(k)}}= \psi^{-1}_{r_{f(k)}}(g-g_1-g_2),  \cr
&  \hskip 4cm \qquad\qquad\qquad\qquad  
\sum_{i=f(k)+1}^{n}\psi_{r_i}(x_{r_i})=g_2, \tau^{*}=k \Bigr \}  \cr}$$
$$\eqalign{
&=q^{-1}\sum_{k=1}^{n^*}\sum_{g_1, g_2 \in G} 
\P\left\{\sum_{i=1}^{f(k)-1}\psi_{r_i}(x_{r_i})=g_1, 
\sum_{i=f(k)+1}^{n}\psi_{r_i}(x_{r_i})=g_2, \tau^{*}=k\right\}  \cr
&=q^{-1}\sum_{k=1}^{n^*} 
\P \left \{ \tau^{*}=k \right \}= 
q^{-1} \P \{\tau^{*}\le n^* \}. \cr
}
$$
Where in the last equalities we have used that $U_{r_{f(k)}}$ is independent
of variables $(x_n:n\neq  r_{f(k)})$ when $\tau^*=k$.
Then,
$$
\P\{S^R_n=g\} =q^{-1}\P\{\tau^{*}\le n^*\} +
\P\{ S^R_n =g, \tau^{*} >n^*\}
$$
and
$$
 \P\{S^R_n=g\} - q^{-1}  = - q ^{-1} \P\{\tau^{*} > n^*\} +
\P\{ S^R_n =g, \tau^{*} >n^*\}.
$$
Using Lemma \equ(4)we get
$$
\left |\P\{S^R_n=g\} - q ^{-1} \right |\le 2 \P \{\tau^{*}>n^*\}\le 
2 \vep(n^*+1). 
$$
\ni b) 
Let $R^J=(R^j:j\in J)$ be a family of sequences,
$(a^{R^j}:j \in J)$ be the family of non-negative sequences, 
$(\psi^{R^j}:j\in J)$ be the corresponding family of 
mappings and $\tilde R^J$ be a family of subsequences verifying 
conditions (H1), (H2), (H3).
Denote by $f_j=f_{\tilde R^j}$ and $\tau_j=\tau^{\tilde R^j}$ for any $j\in J$.
Fix $n \in \N$ and 
put $\tilde n=\tilde n(\tilde R^J)$.  

Take a vector
$\vec k=(k_j:j \in J) \in \{1,...,\tilde n\}^J$.
On the set $\{\tau_j=k_j:j\in J\}$ we define the random
variables 
$$\rho_j(\vec k,n,\underline U)=
\sum_{i \in J} \ \one\{\tilde r_{k_i}^i\in R^j_n\} \ \psi_{f_i(k_i)}(U_{f_i(k_i)}),
\hbox{ for } j \in J.$$
Consider $(g_j':j\in J)\in G^J$. From hypothesis
(H1), (H2), (H3) the system of linear equations
$\rho_j(\vec k,n,\underline U)=g'_j$, $j\in J$, defines  
a system of type (S). Then, by Lemma \equ(402),
there is a unique  $(g_j'':j\in J) \in G^J$ such that
$$\rho_j(\vec k,n,\underline U)=g'_j, \ j\in J\ 
\Leftrightarrow \ U_{f_j(k_j)}=g''_j, \ j \in J. \Eq(403)$$
Let 
$T(\vec k)=(\displaystyle\bigcup_{j\in J} R_n^j)\setminus \{f_j(k_j):j\in J\}.$
It is easy to see that variables $(S_n^{R^j}: j\in J)$ on 
$\{\tau_j=k_j:j\in J\}$ can be written as
$$S_n^{R^j}=\sum_{r \in T(\vec k)\cap R^j_n} \ \psi_r(x_r) \ + \
\rho_j(\vec k,n,\underline U).$$ 
Therefore,
$$\eqalign{
& \P\{S^{R^j}_n=g_j, \tau_j =k_j, \hbox{ for } j\in J\}=\cr 
&\sum_{h_r\in G: \ r\in T(\vec k)} 
\P\{\rho_j(\vec k,n,\underline U)  =g_j- \hskip -0.5 cm \sum_{r\in T(\vec k)\cap R_n^j} 
\psi_r(h_r), 
x_r=h_r, 
\tau_j=k_j, \hbox{ for } j\in J,  
r \in T(\vec k)\} \cr
&=\sum_{h_r \in G:\ r \in T(\vec k)}
\P\{ U_{f_j(k_j)}=g''_j, x_r=h_r, \tau_j=k_j, 
\hbox{ for } \ j \in J,
r \in T(\vec k)\},  \cr
}
$$
where $(g_j'':j \in J) \in G^J$ is given by property \equ(403). 
By independence we conclude that
$$\P\{S_n^{R^j}=g_j, \tau_j=k_j, \hbox{ for } j\in J\}=
q^{-|J|} \P\{\tau_j=k_j, \ j \in J\}.
$$
Hence
$$
\P\{S^{R^j}_n = g_j,\hbox{ for }
j\in J, \ \max\limits_{j\in J} \tau_j\le \tilde n\}=
q^{-|J|}\P\{\max\limits_{j\in J} \tau_j\le \tilde n\},
$$
which together with Lemma  \equ(4) allow us   to deduce that
$$\left |\P\{S^{R^j}_n=g_j, \hbox{ for } j\in J\}-q^{-|J|}\right | \le 2 
\P\{\max\limits_{j\in J}\tau_j>\tilde n\}\le 2|J| \vep(\tilde n +1).$$
\square
\bigskip
Now we can give the proof of the main theorem.
\bigskip
\noindent{\bf Proof of Theorem \equ(1).} 
\bigskip
\bigskip
First, let us introduce some notation. The $p$-expansion of $m\in\N$
is $m=\sum\limits_{i\ge 0} m_ip^i$ with $m_i\in\Z_p$. We denote by
${\cal I}(m)=\{i\in\N :m_i\neq 0\}$ its support
and we denote its elements in decreasing 
order, ${\cal I}(m)=\{\delta_{1,m}>...>\delta_{s_m,m}\}$, where
$s_m=|{\cal I}(m)|$. Now put $m^{(i)}=m_{\delta_{i,m}}$, so $m=\sum\limits^
{s_m}_{i=1} m^{(i)} p^{\delta_{i,m}}$.  Observe that $\delta_{1,m}=$
integer part $(\log m)$, where we take $\log m$ in base $p$.

Since $p$ is a prime number the Lucas' theorem [Lu] asserts that
$$\left[ {m\choose k}\right]_p =\left[ \prod\limits_{i\ge 0}{m_i\choose k_i}
\right]_p, $$
where ${r\choose s}=0$ if $r<s$. In particular
$[{m\choose k}]_p>0$ if and only if $k_i\le m_i$ for all
$i\ge 0$. 

Let us return to the automaton $\varphi$. Since $G$ is Abelian, a simple recurrence
implies
$$
(\varphi^mx)_i=\sum\limits_{k\le m} {m\choose k} \mu^{m-k} \nu ^k x_{k+i}.$$
Observe that
this expression has the form of variables
$S_n^R$ defined before. In this case the mapping
has the shape $ {m \choose k } \mu^{m-k}\nu^k \ g$
which is one-to-one if $\left[ {m \choose k } \right ]_p \neq 0$
since $\mu$ and $\nu$ are relatively primes to $p$.
Then our computations are devoted to show that 
we have enough one-to-one mappings.

In order to make the proof more clear 
we shall first prove that the Ces\`aro mean of the marginal distribution
exists and it is uniform, that means  
$$
\hskip -0.5 cm \pi(g) \doteq\lim\limits_{M\to\infty}{1\over M}\sum\limits^{M-1}_{m=0}
\P_{\underline w}\{(\varphi^mx)_0=g\}\hbox{ \quad exists and verifies \quad}
\pi(g)=q^{-1},\hbox{ for any }g\in G.
$$

Let us fix 
$\alpha\in (0,{1\over 2})$. For
$M>0$ consider the set
$
{\cal R}_M=\{m\le M: |{\cal I}(m)| \ge \alpha\log\log M\}.$ We will
prove that $({\cal R}_M:M\in \N)$ is a sequence of sets of density one,
which means 
${|\{m\le M\}\setminus {\cal R}_M| /  M}
\sobre{\hbox{\rightarrowfill}}{M\to\infty}0$. 
In that purpose we make the 
decomposition \hfill\break
$
\{m\le M\}=\bigcup\limits_{1\le s\le s_M+1} A_{s,m}
$
with 
$$
\eqalign{
&A_{1,M}= \{m\le M: \delta_{1,m}<\delta_{1,M}\},\cr
&A_{s,M}= \{m\le M: \delta_{r,m}=\delta_{r,M}
\hbox{ for } r<s \hbox{ and } \delta_{s,m}<\delta_{s,M}\} \hbox{ for } 
1\le s\le s_M,\cr
&A_{s_M+1,M}= \{M\}. }
$$
Observe that  $|A_{s,M}|=M^{(s)} p^{\delta_{s,M}}$ for $1\le s\le s_M$. 
Take $s^*_M=\sup\{ s:\delta_{s,M}\ge \log\log M\}$. Since
$\delta_{1,M}=$ integer part $(\log M)$, we have $s^*_M\ge 1$. Now,
$$
|\{m\in A_{s,M}: |{\cal I}(m)|\le \alpha\delta_{s,M}\}|\le \sum\limits_{t\le
\alpha\delta_{s,M}} (p-1)^t{\delta_{s,M}\choose t}
$$
$$\le(p-1)^{\alpha\delta_{s,M}}2^{\delta_{s,M}}e^{-2(\alpha-{1\over 2})^2
\delta_{s,M}}.$$
Hence,
$$
|\{m\le M\} \setminus {\cal R}_M|\le\sum\limits_{1\le s\le s^*_M}
(2(p-1)^\alpha)^{\delta_{s,M}}e^{-2(\alpha-{1\over 2})^2
\delta_{s,M}}
+\sum\limits_{s_M^*<s\le s_M} M^{(s)} p^{\delta_{s,M}}+1.
$$
We have
$$
\sum\limits_{s_M^*<s\le s_M} M^{(s)}  p^{\delta_{s,M}}+1\le 
(\log M)^2 +1 .
$$
Take  $\alpha <{p\over 2}(\log(p-1))^{-1}$, then  $p'\doteq 2(p-1)^
\alpha e^{-2(\alpha-{1\over 2})^2} < p$.  Therefore 
$${1 \over M} \sum\limits_{1\le s\le s^*_M}(2(p-1)^\alpha)^{\delta_{s,M}}e^
{-2(\alpha-{1\over 2})^2 \delta_{s,M}} \le {1\over M}
\sum\limits_{1\le s\le s^*_M} p'^{\delta_{s,M}}$$
$$
\le \sum\limits_{1\le s\le s^*_M}\( {p'\over p}\)^{\delta_{s,M}} \le
{p\over p-p'}\({p'\over p}\)^{\log \log M}.
$$

${|\{m\le M\}\setminus {\cal R}_M|\over M}
\sobre{\hbox{\rightarrowfill}}{M\to\infty}0$. 
So $({\cal R}_M:
M\in\N)$ is a sequence of sets of density one. Hence,
$$
\eqalign{\pi(g) &=\lim\limits_{M\to\infty} {1\over M} \sum\limits_{m\in{\cal R}_M}
 {1 \over M} \P_{\underline w} \{ (\varphi^mx)_0 =g\}\cr
&=
\lim\limits_{M\to\infty} {1 \over M} \sum\limits_{m\in {\cal R}_M}
\P_{\underline w} \left\{ \sum\limits_{k\le m} {m \choose k} \mu^{m-k} 
\nu^k x_k=g \right\}.\cr}
$$
 From the Remark,  ${m \choose k}  \ne 0 \ mod \ p$ implies that
the mapping $\psi(g)= {m \choose k} \mu^{m-k}\nu^k\  g $ is one-to-one. Therefore from Lucas' theorem and
Lemma \equ(401) (a)  we get that for any $m \in {\cal R}_M$  
$$\left | \big\{k\le m: {m \choose k} \ mod \ p \ne 0 \big\} \right | 
\ge 2^{\alpha \log\log M}$$
and then 
$$
\left| \P_{\underline w}\left\{\sum\limits_{k\le m} {m \choose k} \mu^{m-k}\nu^k \ x_k=g\right\}-q^{-1}
\right| \le \vep_1(2^{\alpha\log\log M}).$$
Then $\pi(g) = q^{-1}$.
\medskip
Now we are ready to prove the result. Notice that for every 
$(g_j:j<s)\in G^s$ there exists a $(g'_j:j<s)\in G^s$ such that
$$
\{x\in G^\N:(\varphi^nx)_j=g_j\hbox{ for } j<s\}=\{x\in G^\N:
(\varphi^{n+j}x)_0=g'_j\hbox{ for }
j<s\} .
$$
Then it suffices to show that for any finite set $J\subseteq\N$ 
with $0\in J$ and $(g_j:j\in J)\in G^{J}$ it is verified,
$$
\lim\limits_{M\to\infty}{1\over M}\sum\limits_{m\le M}
\P_{\underline w}\{(\varphi^{m+j} x)_0=g_j, j\in J\} =q^{-|J|} .
$$
Introduce the following notation. We put 
$G_m=|\{n\le \delta_{1,m}:m_n<p-1\}|$ and we denote
$$
\{n\le \delta_{1,m}:m_n<p-1\}=\{\beta_{1,m}<\beta_{2,m}<...<\beta_{G_m,m}\}.
$$
Fix $\alpha\in(0,{1\over 2}),\varepsilon\in (0,\alpha),\varepsilon '\in 
(0,{1\over 2}(\alpha-\varepsilon))$. Denote $\ell =\max J$ and define
$$
{\cal R}'_M=\{ m\le M:\log(2(\ell+1))\le G_m\hbox{ and }
\beta_{[\log 2(\ell+1)],m}\le \varepsilon\log\log M\}.
$$
$$
{\cal R}''_M=\{ m\le M:\delta_{1,m}>\varepsilon\log\log M,\
|{\cal I}(m)\cap\{\varepsilon\log\log M \le n \le \delta_{1,m}\}|
\ge \varepsilon ' \log\log M\}.
$$
Both families of sets $({\cal R}'_M:M\in\N)$, $({\cal R}''_M:M\in \N)$ are 
of density 1.
\smallskip
Now for any family of sets $(\tilde {\cal R}_M:M\in \N)$ with 
$\tilde {\cal R}_M\subseteq \{m\le M\}$, we put $\tilde
{\cal R}_{M,J}=\{m\le M:
m+j\in \tilde{\cal R}_M$ for $j\in J\}$. If $(\tilde {\cal R}_M:M\in\N)$ 
is of density 1 then also
$(\tilde {\cal R}_{M,J}:M\in \N)$ is of density 1. Hence $({\cal R}_{M,J}:M\in \N)$, 
$({\cal R}'_{M,J}:M\in \N)$,
$({\cal R}''_{M,J}:M\in \N)$ are sequences of density 1.
\smallskip

Let $m\in {\cal R}'_{M,J}\cap {\cal R}''_{M,J}$. We denote 
${\cal I}_+(m+j)={\cal I}(m+j)\cap\{n>
\varepsilon \log\log M\}$, and ${\cal I}_-(m+j)={\cal I}(m+j)
\cap\{n\le \varepsilon\log\log M\}$. From the definition of ${\cal R}'_M$
we have that ${\cal I}_+(m+j)={\cal I}_+(m)$ for $j\in J$. Put
${\cal C}_+(m+j)=\{(m+j)_i:i\in {\cal I}_+(m+j)\}$ and
${\cal C}_-(m+j)=\{(m+j)_i:i\in {\cal I}_-(m+j)\}$ for $j\in J$. We have
${\cal C}_+(m+j)={\cal C}_+(m)$ for $j\in J$, and the sets
$({\cal C}_-(m+j):j\in J)$ are all different between them. Define for $j\in J$
$$
\tilde {\cal R}^j \! = \!
\{ 
k \le m+j: {\cal I}(k) \! \subseteq \! {\cal I}(m+j), k_i \! \le \! m_i
\hbox{ for } i \in \! {\cal I}_+(m),
n_i \! =  \! (m+j)_i \hbox{ for } i \in  \! {\cal I}_-(m+j) \}.
$$
The family $(\tilde {\cal R}^j:j\in J)$ is disjoint because the sets  
$({\cal C}_-(m+j):j\in J)$ are different. Moreover 
$|\tilde {\cal R}^j|\ge 2^{\varepsilon' \log\log M}$. 
\medskip
 From Lemma \equ(401) (b) and the Remark 
we get the result. In fact for every $m\in {\cal R}'_{M,J}\cap {\cal R}''_{M,J}$
we have that
$$
(\varphi^{m+j}x)_0=\sum\limits^{m+j}_{k=0} {m+j \choose k} \mu^{m+j-k}\nu^k x_k
$$
and the sequences 
$(\tilde {\cal R}^j:j\in J)$ satisfies conditions (H1),(H2),(H3). Indeed,
property (H1) follows from  
$\tilde {\cal R}^j \subset \{k\le m+j :{m+j \choose k} \ mod \ p   >0 \}
$,
they are disjoint, and if $k \in \tilde R^j$ then 
${m+j' \choose k} \ mod \ p=0$ for every $j'<j$ in  $J$ which shows property
(H3).
Then, from Lemma \equ(401) (b),  for any such $m$
$$
\left |\P_{\underline w}\{x:(\varphi^{m+j}x)_0=g_j, j\in J\}-q^{-|J|}\right |
\le \vep_J( 
2^{\varepsilon' \log\log M}).
$$
Then the theorem is shown. \square

\bigskip
\ni{\bf Acknowledgments.} Alejandro Maass and Servet Mart\'inez acknowledge 
financial support
 from C\'atedra Presidencial fellowship and Fondecyt grants 1980657 and
1970506.
Pablo A. Ferrari is partially supported by FAPESP
(Projeto Tem\'atico), CNPq (Bolsa de aux\'\i lio \`a pesquisa) and
FINEP (Projeto N\'ucleos de Excel\^encia).

\bigskip
\vfill\eject
\ni {\bf References.}
\bigskip
\item{[AN]} K.B. Athreya, P. Ney, {\it A new approach to the limit theory
of recurrent Markov chains}, Transactions of the AMS 248, 493--501 (1978).
\medskip

\item{[B]} H. Berbee (1987)  {\it Chains with infinite connections: Uniqueness
and Markov representations}, 
Probab. Theory Related Fields 76 (1987), no. 2, 243--253.  
\medskip

\medskip

\item{[BFG]} X. Bressaud, R. Fern\'andez, A. Galves (1998) {\it Decay of
correlations for non H\"olderian dynamics. A coupling approach. Preprint.

\item{[FG]} P. Ferrari, A. Galves (1997) 
{\it Acoplamento em processos estocásticos}, 
21 Colóquio Brasileiro de Matemática [21th Brazilian Mathematics Colloquium],
IMPA, Rio de Janeiro. Available in
http://www.ime.usp.br/~pablo/abstracts/libro.html. 
\medskip

\item{[FMM]} P. Ferrari, A. Maass, S. Mart\'inez (1998)
{\it Cesaro mean distributionn of group automata 
starting from Markov measures}, 
Unpublished.
\medskip

\item{[L]} D. Lind, {\it Applications of ergodic theory and sofic systems 
to cellular automata}, Physica D 10, 36-44 (1984).
\medskip

\item{[Lu]} E. Lucas, {\it Sur les congruences des nombres eul\'eriens et des
coefficients diff\'erentiels des fonctions trigonom\'etriques,
suivant un module premier}, Bulletin de la Soc. Math\'ematique 
de France 6, 49--54 (1878).

\medskip

\item{[MM]} A. Maass, S. Mart\'inez, {\it On Ces\`aro limit distribution of
a class of permutative cellular automata}, Journal of Statistical Physics
90, 435--452 (1998).  

\item{[NN]} Ney, P.; Nummelin, E. {\it Regeneration for chains with infinite
memory}, Probab. Theory Related Fields 96 (1993), no. 4, 503--520.

\bye